\documentclass[11pt,leqno]{article}

\usepackage{a4wide,amssymb,amsmath,amscd,amsthm,amsfonts,graphicx}

\usepackage[dvips]{pstcol}

\usepackage[all]{xy}
\newcommand{\Magma}{\textsc{Magma}}

\newtheorem{theorem}{Theorem}
\newtheorem{lemma}[theorem]{Lemma}
\newtheorem{proposition}[theorem]{Proposition}
\newtheorem{corollary}[theorem]{Corollary}
\newtheorem{remark}[theorem]{Remark}

\def\semi{\hbox{ $\times $ \kern-.972em \raise.12719em\hboX{ $_{^|}$}  }}

\def\be{\begin{enumerate}}
\def\ee{\end{enumerate}}
\def\bi{\begin{itemize}}
\def\ei{\end{itemize}}

\def\pbraid{\Sigma}

\title{Conjugacy in Garside Groups III: Periodic braids
  }
\author{Joan S. Birman\footnote{Partially
supported by the U.S.National Science Foundation, under Grant
DMS-0405586.} \and  Volker Gebhardt \and Juan
Gonz\'alez-Meneses\footnote{Partially supported by
MTM2004-07203-C02-01 and FEDER.}}
\date{February 19, 2007}

\begin{document}

\maketitle

\begin{abstract}
An element in Artin's braid group $B_n$ is said to be periodic if some power
of it lies in the center of $B_n$. In this paper we prove that all
previously known algorithms for solving the conjugacy search problem in $B_n$
are exponential in the braid index $n$ for the special case of periodic
braids.
We overcome this difficulty by putting to work several
known isomorphisms between Garside structures in the braid group $B_n$ and
other Garside groups.  This allows us to obtain a polynomial solution to the
original problem in the spirit of the previously known algorithms.

This paper is the third in a series of papers by the same authors about the
conjugacy problem in Garside groups. They have a unified goal: the development
of a polynomial algorithm for the conjugacy decision and search problems in
$B_n$, which generalizes to other Garside groups whenever possible. It is our
hope that the methods introduced here will allow the generalization of the
results in this paper to all Artin-Tits groups of spherical type.
\end{abstract}

\section{Introduction}
Given a group, a solution to the {\it conjugacy decision problem} is an
algorithm that determines whether two given elements are conjugate
or not. On the other hand, a solution to the {\it conjugacy search
  problem} is an algorithm that finds a conjugating element for a
given pair of conjugate elements. In $\S$1.4 of~\cite{BGGM-I} we presented a
project to find a polynomial solution to the conjugacy decision
problem and the conjugacy search problem in the particular case of
Artin's braid group, that is, the Artin-Tits group of type
${\bf A}_{n-1}$, with its classical  or {\it Artin} presentation~\cite{Artin}:
\begin{equation}
\label{E:classical presentation} B_n^{A} :  \left< \sigma_1,\ldots,
\sigma_{n-1} \left|
\begin{array}{ll}
 \sigma_i \sigma_j = \sigma_j \sigma_i  & \mbox{if } \; |i-j|>1,  \\
 \sigma_i \sigma_j \sigma_i = \sigma_j \sigma_i \sigma_j &  \mbox{if }\; |i-j|=1.\ \end{array}\right.
 \right>.
\end{equation}
One of the steps in the mentioned project asks for a polynomial
solution to the above conjugacy problems for special type of
elements in the braid groups, called {\it periodic braids}. This is achieved in
the present paper. More precisely, if we denote by $|w|$ the letter
length of a word $w$ in $\sigma_1,\ldots,\sigma_{n-1}$ and their
inverses, we will prove:

\begin{theorem}  \label{T:main theorem}
Let $w_X$ and $w_Y$ be two words in the generators
$\sigma_1,\ldots,\sigma_{n-1}$ and their inverses, representing two
braids $X,Y \in B_n^A$, and let $l=\max\{|w_X|,|w_Y|\}$. Then there is
an algorithm of complexity $O(l^3n^2\log n)$ which does the
following. \be
\item [{\rm (1)}] It determines whether $X$ and $Y$ are periodic.
\item [{\rm (2)}] If yes, it determines whether they are conjugate.
\item [{\rm (3)}] If yes, it finds a braid $C\in B_n^A$ such that $Y = C^{-1}XC$.
\ee
\end{theorem}

Here is a guide to this paper. In Section~\ref{S:Definitions and known results}, we will review what is known and explain why
steps $(1)$ and $(2)$ of Theorem~\ref{T:main theorem}  follow easily from the work
in~\cite{Eilenberg,Ker,GM}.   On the other hand, in Section~\ref{S:Known algorithms} we show that the previously
known solutions to the conjugacy search problem in the Artin-Tits group of type
${\bf A}_{n-1}$  present unexpected difficultites, which result in exponential complexity for periodic braids. Thus they do not meet the requirements of Theorem 1.  

A new idea allows us to overcome the difficulty.  We have shown that the approach using the classical Garside structure does not work. The new idea is to put to work the other known Garside structure on the braid groups and in addition to consider a certain subgroup of the braid group that arises in the course of our work, and use two known Garside structures on it.
This is accomplished in Section~\ref{S:proof}, where we give a solution to the conjugacy search
problem for periodic braids which has the stated polynomial
complexity.   Section~\ref{S:proof} divides naturally into two subsections, according to whether a given periodic braid is conjugate to a power of $\delta$ or $\varepsilon$,  two braids that are defined in 
Section~\ref{S:Definitions and known results} below. The proof in the two cases are treated in Sections \ref{SS:powers of delta} and \ref{SS:powers of varepsilon} respectively.
Finally, in Section~\ref{S:Timing results} we compare
actual running times of the algorithms developed in
Section~\ref{S:proof} to the ones of the best previously known
algorithm.

\begin{remark} {\rm We learned from D. Bessis that he has
characterized the conjugacy classes of periodic elements for all
Artin-Tits groups of spherical type. We hope that this
characterization will allow the generalization of both the techniques and the results of this paper
to all other Artin-Tits groups of spherical type.}
\end{remark}

\noindent {\bf Acknowledgements:} We are grateful to D. Bessis
for useful discussions about his work in~\cite{Bessis} and his
forthcoming results, to J. Michel for pointing out that our
Corollaries~\ref{C:size of USS(delta)} and \ref{C:size of
USS(epsilon)} were known to specialists in Coxeter groups, and
also to H. Morton for showing us the algorithm in~\cite{M-H}.

\section{Known results imply steps (1) and (2) of Theorem~\ref{T:main theorem}}
\label{S:Definitions and known results}

Our work begins with a review of known results. Garside groups were introduced by Dehornoy and Paris in~\cite{D-P}.
The main examples of Garside groups are Artin-Tits groups of
spherical type, in particular, Artin braid groups.   In this paper we will use
two known Garside structures in the Artin-Tits group of type
${\bf A}_{n-1}$, and also one
Garside structure in the Artin-Tits group of type ${\bf B}_m$.

Although we refer to~\cite{BGGM-I} for a detailed description of
Garside structures, we recall here that such a structure in a group
$G$ is given by a lattice order on its elements, together with a
distinguished element of $G$, called the {\it Garside element},
which is usually denoted by $\Delta$. This partial order and this
element $\Delta$ must satisfy several suitable
conditions~\cite{BGGM-I}.

The classical Garside structure in the braid groups is related to
the presentation~(\ref{E:classical presentation}). The {\it
positive} braids are those which can be written as a word in
$\sigma_1,\ldots,\sigma_{n-1}$ (not using their inverses). The
lattice order is defined by saying that $X\preccurlyeq Y$ if
$X^{-1}Y$ is a positive braid (we will say that $X$ is a {\it
prefix} of $Y$). There are special elements called {\it simple
braids} which are those positive braids in which any two strands
cross {\it at most} once. The Garside element $\Delta$ is the
positive braid in which any two strands cross {\it exactly} once,
that is, $\Delta=\sigma_1
(\sigma_2\sigma_1)(\sigma_3\sigma_2\sigma_1)\cdots(\sigma_{n-1}\cdots
\sigma_1)$. It is also called the {\it half twist}, since its
geometrical representation corresponds to a half twist of the $n$
strands. For every braid $X\in B_n^A$, given as a word of letter
length $l$, there exists a {\it left normal form}, which is a unique
way to decompose the braid as $X=\Delta^p x_1\cdots x_r$, where $p$
is maximal and each $x_i$ is a simple braid, namely the maximal
simple prefix of $x_i\cdots x_r$. This left normal form can be
computed in time $O(l^2n\log n)$~\cite{Epstein}.

Artin proved  in \cite{Artin} that the center of $B_n^A$ is infinite
cyclic and generated by the {\it full twist} $\Delta^2 =
(\sigma_1\sigma_2\cdots\sigma_{n-1})^n$ of the braid strands.  If
the braid group is regarded as the mapping class group of the
$n$-times punctured disc $\mathbb D^2_n$, then $\Delta^2$ is a Dehn
twist about a curve which lies in a collar neighborhood of the
boundary $\partial\mathbb D^2_n$ and is parallel to it.  An element
$X\in B_n^A$  is said to be {\it periodic} if some power of $X$ is a
power of $\Delta^2$.

\begin{figure}[ht]
\centerline{\includegraphics{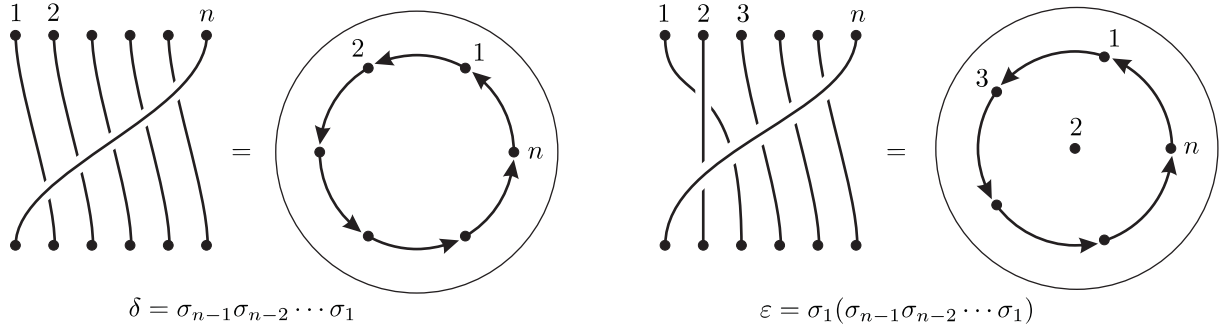}} \caption{The
periodic elements $\delta$ and $\varepsilon$.}
\label{F:delta_epsilon}
\end{figure}

Periodic braids can be thought of as {\it rotations} of the disc.
Indeed, there is a classical result by Eilenberg~\cite{Eilenberg}
and K\'er\'ekj\'art\'o\cite{Ker} (see also~\cite{C-K}) showing that
an automorphism of the disc which is a root of the identity (a
periodic automorphism) is conjugate to a rotation. Since a finite
order mapping class can always be realized by a finite order
homeomorphism~\cite{Kerck}, this implies that a periodic braid is
conjugate to a rotation. It is not difficult to see that a braid can
be represented by a rotation of $\mathbb D^2$ if and only if it is
conjugate to a power of one of the two braids represented in
Figure~\ref{F:delta_epsilon}, that is,
$\delta=\sigma_{n-1}\sigma_{n-2}\cdots \sigma_1$ and
$\varepsilon=\sigma_1 (\sigma_{n-1}\sigma_{n-2}\cdots \sigma_1)$.
(If we need to specify the number of strands, we will write
$\delta=\delta_n$ and $\varepsilon=\varepsilon_n$.) 

\begin{remark} {\rm   The braid $\varepsilon$ defined in Figure~\ref{F:delta_epsilon}
 has a fixed strand,
namely strand 2. There are, to be sure, braids which are conjugate to
$\varepsilon$ in which the fixed strand is the first one or the last
one, seemingly more natural choices. However, $\varepsilon$ is a simple braid, and 
(as we shall prove in Proposition~\ref{P:USS of epsilon}  below) there is no simple braid which is 
conjugate to $\varepsilon$ and which fixes either the first or the last
strand. This is why we decided to use $\varepsilon$, which fixes the
second strand, as a representative of its conjugacy class.  And this
is also the reason why, in Section~\ref{SS:powers of varepsilon}  below, we identify the
Artin-Tits group of type $\mathbf B_{n-1}$ with the subgroup of the $n$-strand braid group
formed by those braids which fix the {\it second} strand, a choice
that will surely seem awkward to specialists.
} \end{remark}

The theorem of
Eilenberg and K\'er\'ekj\'art\'o can then be restated as follows.

\begin{theorem}{\rm \cite{Eilenberg,Ker}}\label{T:Eilenberg-Kerekjarto}
A braid $X$ is periodic if and only if it is conjugate to a power of either
$\delta$ or $\varepsilon$.
\end{theorem}

\medskip
Notice that $\delta^n=\varepsilon^{n-1}=\Delta^2$. Since $\Delta^2$
belongs to the center of $B_n^A$, this immediately gives an efficient
algorithm to check whether a braid is periodic.

\begin{corollary}\label{C:Eilenberg-Kerekjarto}
A braid $X\in B_n^A$ is periodic if and only if either $X^{n-1}$ or
$X^n$ is a power of $\Delta^2$.
\end{corollary}

\begin{proof}
We only need to prove that the condition is necessary. Suppose that
$X$ is periodic. By Theorem~\ref{T:Eilenberg-Kerekjarto}, $X$ is
conjugate to a power of either $\delta$ or $\varepsilon$. In the
first case, $X=C^{-1}\delta^k C$ for some $C\in B_n^A$. Then
$X^n=C^{-1}\delta^{kn}C = C^{-1}\Delta^{2k}C = \Delta^{2k}$, where
the last equality holds since $\Delta^2$ is central. In the second
case, $X=C^{-1}\varepsilon^k C$, so that $X^{n-1}= C^{-1}
\varepsilon^{k(n-1)}C = C^{-1} \Delta^{2k} C = \Delta^{2k}$.
\end{proof}

After this result, one can determine whether $X$ is periodic, and
also find the power of $\delta$ or $\varepsilon$ which is conjugate
to $X$, by the following algorithm.

\medskip
\noindent {\bf Algorithm A.}

Input: A word $w$ in Artin generators and their inverses
representing a braid $X\in B_n^A$.

\begin{enumerate}

  \item Compute the left normal form of $X^{n-1}$.

  \item If it is equal to $\Delta^{2k}$, return {\it `$X$ is periodic and conjugate to $\varepsilon^k$'.}

  \item Compute the left normal form of $X^n$.

  \item If it is equal to $\Delta^{2k}$, return {\it `$X$ is periodic and conjugate to $\delta^k$'.}

  \item Return {\it `$X$ is not periodic'.}

\end{enumerate}

\medskip

\begin{proposition}~\label{P:Algorithm A}
The complexity of Algorithm A is $O(l^2n^3\log n)$, where $l$ is the
letter length of $w$.
\end{proposition}

\begin{proof}
Algorithm A computes two normal forms of words whose lengths are at
most $nl$. By~\cite{Epstein}, these computations have complexity
$O((nl)^2 n\log n)$, and the result follows.
\end{proof}

We remark that if one knows, a priori, that the braid $X$ is
periodic, then one can determine the power of $\delta$ or
$\varepsilon$ which is conjugate to $X$ by a faster method: Observe
that the exponent sum of a braid $X$, written as a word in the
generators $\sigma_1,\ldots,\sigma_{n-1}$ and their inverses is well
defined, since the relations in~(\ref{E:classical presentation}) are
homogeneous. The exponent sum is furthermore invariant under
conjugacy, hence every conjugate of $\delta^k$ has exponent sum
$k(n-1)$, whereas every conjugate of $\varepsilon^k$ has exponent
sum $kn$. Moreover, the exponent sum determines the conjugacy class
of a periodic braid:

\begin{lemma} \label{L:exponent sum} {\rm (Proposition 4.2 of \cite{GM})}
Let $X$ be a periodic braid. Then $X$ is conjugate to $\delta^k$
(resp. $\varepsilon^k$) if and only if $X$ has exponent sum $k(n-1)$
(resp. $kn$).
\end{lemma}

Computing the exponent sum of a word of length $l$ has complexity
$O(l)$. Hence, once it is known that two given braids are periodic,
the conjugacy decision problem takes linear time.

\section{Known algorithms are not efficient for periodic braids}
\label{S:Known algorithms}

We have already determined all conjugacy classes of periodic braids,
and we have seen that the conjugacy decision problem for these
braids can be solved very fast. It is then natural to wonder whether
this is also true for the conjugacy search problem. The first
natural question is: Are the existing algorithms for the conjugacy
search problem efficient for periodic braids?

The best known algorithm to solve the conjugacy decision problem and
also the conjugacy search problem in braid groups (and in every
Garside group) is the one in~\cite{Gebhardt}, which consists of
computing the {\it ultra summit set} of a braid, defined as follows.
Denote by $\tau$ the inner automorphism that is defined by
conjugation by $\Delta$. Given $Y\in B_n^A$ whose left normal form is
$\Delta^p y_1\cdots y_r$,
we define its {\it canonical length} as $\ell(Y)=r$,
and call the conjugates
$\mathbf c(Y)=\Delta^p y_2 \cdots y_r \tau^{-p}(y_1)$ and
$\mathbf d(Y)=\Delta^p \tau^p(y_r) y_2 \cdots y_{r-1}$ of $Y$
its {\it cycling} respectively its {\it decycling}.
For every $X\in B_n^A$,
the ultra summit set $USS(X)$ is the set of conjugates $Y$ of $X$
such that $\ell(Y)$ is minimal and $\mathbf c^t(Y)=Y$ for some
$t\geq 1$.  It is explained in~\cite{Gebhardt} how the computation
of $USS(X)$ solves the conjugacy decision and search problems in
Garside groups.

The complexity of the conjugacy search algorithm given
in~\cite{Gebhardt} is proportional to the size of $USS(X)$, so if
one is interested in complexity, it is essential to know how large
the ultra summit sets of periodic braids are. If they turned out to
be small, the algorithm in~\cite{Gebhardt} would be efficient, but
we will see in this section that the sizes of ultra summit sets of
periodic braids are in general exponential in $n$.

More precisely, it was shown by Coxeter in 1934~\cite[Theorem
11]{Coxeter}, that in any finite Coxeter group, any two elements
which are the product of all standard generators, in arbitrary
order, are conjugate. Applied to our case, one sees that the
elements of $USS(\delta)$ are in bijection with the elements of the
above kind, in the symmetric group $\Sigma_n$. One can count the
number of different elements, and it follows that $\#(USS(\delta))=
2^{n-2}$. The same result is shown in~\cite[Chapter V, \S6.
Proposition 1]{Bourbaki}, in the more general case in which the
Coxeter group is defined by a tree, and also in \cite[Lemma
3.2]{Shi} and in~\cite[Theorem 2]{M-H}. Moreover, it can be seen
from the proof in~\cite{Bourbaki} that any two elements in
$USS(\delta)$ are conjugate by a sequence of special conjugations,
that we denote {\it partial cyclings} in~\cite{BGGM-I}.

Concerning the elements in $USS(\varepsilon)$,
in~\cite[Proposition 9.1]{Digne-Michel} it is shown that any
two such elements are conjugate by a sequence of partial
cyclings. It also follows from~\cite{Digne-Michel} that every
element in $USS(\varepsilon)$ is represented by a word of
length $n$, which is the product of all $n-1$ generators, in
some order, with one of the generators repeated. One can also
count the number of different elements of this kind, to obtain
that $\#(USS(\varepsilon))= (n-2)2^{n-3}$.

The above arguments show that the sizes of $USS(\delta)$ and
$USS(\varepsilon)$ are exponential with respect to the number
of strands, hence the algorithm in~\cite{Gebhardt} is not
polynomial for conjugates of these braids. In this paper we
shall study $USS(\delta)$ and $USS(\varepsilon)$ in a new way.
More precisely, in Corollaries~\ref{C:size of USS(delta)} and
\ref{C:size of USS(epsilon)} we will show that
$\#(USS(\delta))= 2^{n-2}$ and $\#(USS(\varepsilon))=
(n-2)2^{n-3}$ just by looking at the permutations induced by
their elements. This will also provide a fast solution to the
conjugacy search problem in the particular cases of conjugates
of $\delta$ or $\varepsilon$.

Once shown that the algorithm in~\cite{Gebhardt} is not polynomial,
in general, for periodic braids, in Section~\ref{S:proof} we will
give a procedure to solve the conjugacy search problem for all
periodic braids in polynomial time.

Let us then study the ultra summit sets of $\delta$ and
$\varepsilon$. First, we recall that the factors in a left normal
form are {\it simple braids}, which are in bijection with the
elements of the symmetric group $\Sigma_n$. More precisely, every
braid $X$, being a mapping class group of the $n$-times punctured
disc, determines a permutation $\pi_X$ of the $n$ punctures.
Conversely, there is exactly one simple braid for each permutation.
We will then determine simple elements by their permutations,
written as a product of disjoint cycles. For instance, the
permutation associated to $\delta$ is $\pi_{\delta}=(1 \ 2\ \cdots\
n)$, and the permutation associated to $\varepsilon$ is
$\pi_{\varepsilon}= (2)(1\ 3\ 4 \ \cdots \ n )$.

\begin{remark} {\rm  Although we described braids as mapping
classes, we will not adopt the usual convention for compositions of
maps. We consider braids as acting on the punctures from the right.
This means that the braid $\sigma_1\sigma_2$ first swaps the
punctures in positions $1$ and $2$, and then the punctures in
positions $2$ and $3$. Hence $\pi_{\sigma_1\sigma_2}=(1 3 2)$.}
\end{remark}

\begin{remark} {\rm  The permutation associated to a simple braid
$s$ determines the pairs of strands that cross in $s$. More
precisely, two strands $i$ and $j$ ($i<j$) cross in $s$ if and only
if the induced permutation reverses their order, that is, if
$\pi_s(i)>\pi_s(j)$.}
\end{remark}

For simplicity of notation let us define, for $1\leq i<j \leq n$,
the braids $\sigma_{[i\rightarrow j]}=\sigma_i\sigma_{i+1}\cdots
\sigma_{j-1}$ and $\sigma_{[j\rightarrow
i]}=\sigma_{j-1}\sigma_{j-2}\cdots \sigma_i$. Notice that
$\sigma_{[k\rightarrow l]}$ (no matter which subindex is bigger) is
the shortest positive braid sending the puncture $k$ to the puncture
$l$.

Let us characterize the elements in $USS(\delta)$.

\begin{proposition}\label{P:USS of delta}
An element $s\in B_n^A$ belongs to $USS(\delta)$ if and only if it is
simple and its permutation $\pi_s$ is a cycle of the form:
$$
 \pi_s= (1 \; u_1 \; u_2 \cdots u_r \; n \; d_t \; d_{t-1} \cdots d_1),
$$
for some $u_1<u_2<\cdots < u_r$ and some $d_t>d_{t-1}>\cdots >d_1$,
with $r,t\geq 0$ and $r+t+2=n$. Moreover, in this case $\alpha^{-1}
s \alpha = \delta$, where
$$
\alpha= \sigma_{[d_1\rightarrow 1]}\:\sigma_{[d_2\rightarrow
1]}\cdots \sigma_{[d_t\rightarrow 1]}.
$$
\end{proposition}

\begin{figure}[ht]
\centerline{\includegraphics{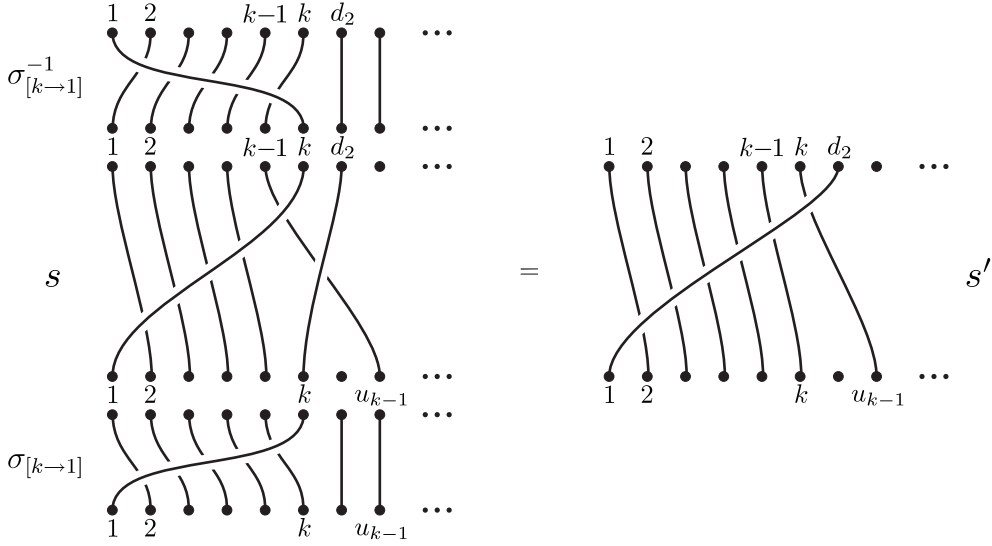}}
\caption{Conjugating $s$ to $s'$.} \label{F:conjugate_of_delta}
\end{figure}

\begin{proof}
First notice that, since $\delta$ is simple, all elements in
$USS(\delta)$ are simple, so that by the definition of a simple
element they can be characterized by their permutations. Actually,
$USS(\delta)$ is the set of simple conjugates of $\delta$.  Notice
also that $\pi_{\delta}$ is a single cycle of length $n$. Since
conjugation of braids in $B_n^A$ implies conjugation of their
corresponding permutations, it follows that the elements in
$USS(\delta)$, which are conjugates of $\delta$, are simple elements
determined by a cycle of length $n$. Moreover, if $s\in USS(\delta)$
then $s^n=\Delta^2$, which is a positive braid in which any two
strands cross exactly twice.

Let $s\in USS(\delta)$. Its permutation can be written as $\pi_s= (1
\; u_1 \; u_2 \ \cdots \ u_r \; n \; d_t \; d_{t-1} \ \cdots\ d_1)$,
where $r,t\geq 0$ and $r+t+2=n$. We must show that $u_1<\cdots <
u_r$ and $d_t>\cdots>d_1$. See in
Figure~\ref{F:conj_and_notconj_of_delta} an example of two simple
braids whose permutations are cycles of length $n$, so the
permutations are conjugate in the symmetric group, but one of the
braids satisfies the above inequalities and the other one does not.

\begin{figure}[ht]
\centerline{\includegraphics{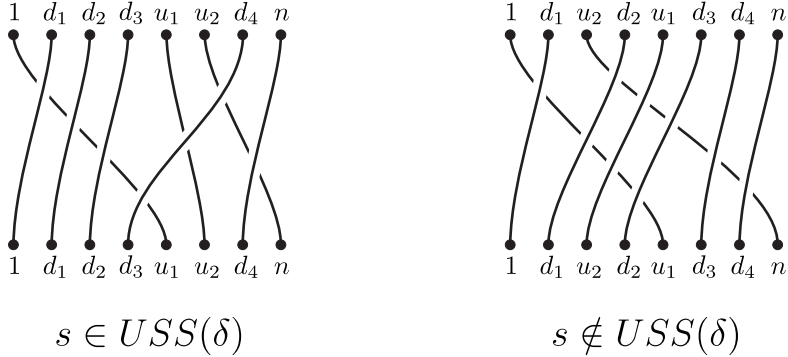}}
\caption{Two simple braids in $B_8^A$ whose permutations are cycles of
length $8$. By Proposition~\ref{P:USS of delta}, the first one is
conjugate to $\delta$ and the second one is not. Notice that the
exponent sum of the second one (i.e. the number of crossings or the
letter length, in this case) is 9, while the exponent sum of
conjugates of $\delta\in B_8^A$ is 7.}
\label{F:conj_and_notconj_of_delta}
\end{figure}

Suppose that $u_i>u_{i+1}$ for some $i$, where $1\le i<r$, and
consider the strands $1$ and $u_{1}$.  We will see that these two
strands cross more than twice in $s^n$. Indeed, one has $1<u_{1}$,
but in $s^i$ these strands end at $u_i$ and $u_{i+1}$, respectively.
Since $u_i>u_{i+1}$, this means that they have crossed at least once
in $s^i$. Now in $s^r$ these two strands end at $u_r$ and $n$,
respectively, and since $u_r$ is necessarily less than $n$, they
have crossed again. Next, in $s^{r+1}$ they end at $n$ and $d_t$ (or
$n$ and $1$ if there are no $d_j$'s), so they have crossed one more
time. This means that in $s^{r+1}$ the strands $1$ and $u_1$ cross
at least three times, showing that $s^n$ cannot be equal to
$\Delta^2$, a contradiction. Therefore $u_1<\cdots<u_r$. Similarly,
if we had $d_{i+1}<d_i$ for some $i$, then strands $n$ and $d_t$
would cross more than twice in $s^n$, which is impossible. Therefore
$d_t>\cdots>d_1$.

Conversely, suppose that $s$ is simple and $\pi_s= (1 \; u_1 \; u_2
\ \cdots\  u_r \; n \; d_t \; d_{t-1}\  \cdots\  d_1)$ for some
$u_1<\cdots <u_r$ and $d_t>\cdots >d_1$. We will show that $s$ is
conjugate to $\delta$ in a constructive way, by finding a
conjugating element. First notice that if $t=0$ then $\pi_s=(1\ 2 \
\cdots \ n) =  \pi_\delta$. Since simple elements are determined by
their permutations, this means that $s=\delta$. Hence we can assume
that $t>0$. Denote $k=d_1$. One has
$$
\pi_s=(\underline {1\ 2\ \cdots \ k-1} \ u_{k-1} \ \cdots \ u_r \ n
\ d_t\ \cdots \ d_2 \underline{\ k}\: ).
$$
A schematic picture of the first $k$ strands of $s$ can be seen in
Figure~\ref{F:conjugate_of_delta}. We will conjugate $s$ by
$\sigma_{[k\rightarrow 1]}$, so we consider $s'=
\sigma_{[k\rightarrow 1]}^{-1}\: s \:\sigma_{[k\rightarrow 1]}$.
Recall that two strands $i$ and $j$ ($i<j$) cross in $s$ if and only
if $\pi_s(i)>\pi_s(j)$. Then we can easily check that the strand of
$s$ ending at $k$ (that is, the strand $d_{2}$ if $t>1$ or the
strand $n$ if $t=1$) does not cross the strands ending at
$1,2,\ldots,k-1$ (that is, the strands $k,1,2,\ldots,k-2$,
respectively). This implies that $s\:\sigma_{[k\rightarrow 1]}$ is a
simple braid. Moreover, one can also check that the strand $k$ of
$s$ (thus the strand $k$ of $s \:\sigma_{[k\rightarrow 1]}$) crosses
the strands $k-1,k-2,\ldots,1$, hence $s'= \sigma_{[k\rightarrow
1]}^{-1}\: s \:\sigma_{[k\rightarrow 1]}$ is a simple braid.  Since
the permutation associated to $\sigma_{[k\rightarrow 1]}$ is $(1\ 2
\ \cdots \ k)$, it follows that
$$
\pi_{s'}=(\underline {1\ 2\ \cdots \ k-1 \ k} \ u_{k-1} \ \cdots \
u_r \ n \ d_t\ \cdots \ d_{2}).
$$
We can continue this process, by recurrence on $t$, conjugating by
elements of the form $\sigma_{[d_i\rightarrow 1]}$ and obtaining new
simple conjugates of $s$ whose permutations have more indices
between $1$ and $n$ at each step, until we get the permutation $(1\
2 \ \cdots\ n)$, that is, until we obtain $\delta$. In this way we
have shown that if $s$ is a simple element with the permutation
given in the statement, then $\alpha^{-1} s \alpha = \delta$, where
$$
\alpha= \sigma_{[d_1\rightarrow 1]}\:\sigma_{[d_2\rightarrow
1]}\cdots \sigma_{[d_t\rightarrow 1]}.
$$
Therefore, we have determined the elements in $USS(\delta)$ in terms
of their permutations.
\end{proof}

\begin{remark} {\rm  The above element $\alpha$ is simple, hence
all elements in $USS(\delta)$ are conjugate to $\delta$ by a simple
element.}
\end{remark}

\begin{corollary}\label{C:size of USS(delta)}
If $\delta=\sigma_{n-1}\cdots\sigma_1\in B_n^A$ then $\#
(USS(\delta))= 2^{n-2}$.
\end{corollary}

\begin{proof}
The elements in $USS(\delta)$ are characterized by the permutation
given in the above result, which is itself characterized by the
sequence $1<u_1<\cdots <u_r<n$. The number of possible sequences is
equal to the number of subsets of $\{2,\ldots,n-1\}$ which is
precisely $2^{n-2}$.
\end{proof}

Now let us do the same for $USS(\varepsilon)$.

\begin{proposition}\label{P:USS of epsilon}
An element $s\in B_n^A$ belongs to $USS(\varepsilon)$ if and only if
it is simple and
$$
 \pi_s= (a)(1 \ u_1 \ u_2 \ \cdots\  u_r \ n \ d_t \ d_{t-1} \ \cdots \ d_1),
$$
for some $u_1<u_2<\cdots < u_r$ and some $d_t>d_{t-1}>\cdots >d_1$,
with $r,t\geq 0$ and $r+t+3=n$. Notice that $a\neq 1,n$.  Moreover,
in this case one has $\beta^{-1} s \beta =\varepsilon$, where
$$
 \beta =\sigma_{[d_1\rightarrow 1]}\:\sigma_{[d_2\rightarrow
1]}\cdots \sigma_{[d_t\rightarrow 1]}\: \sigma_{[b\rightarrow 2]}
$$
and $b=a+t-\max\{i\,:\,d_i<a \}$,
\end{proposition}

\begin{proof}
Since $\varepsilon$ is simple, the elements of $USS(\varepsilon)$
are precisely the simple conjugates of $\varepsilon$; in particular,
$USS(\varepsilon)$ consists of simple elements whose permutation is
the product of a cycle of length 1 (a fixed point) and a cycle of
length $n-1$. Moreover, if $s\in USS(\varepsilon)$ then
$s^{n-1}=\Delta^2$, where any two strands cross exactly twice.

Let $s\in USS(\varepsilon)$, and let $ \pi_s= (a)(x_1  \ \cdots\
x_{n-1})$. If $a=1$ then the first strand of $s$ does not cross any
other strand. This means that we can write $s$ as a word in Artin
generators in which the letter $\sigma_1$ does not appear. But in
that case every power of $s$ would satisfy the same property. In
particular, the first strand of $s^{n-1}=\Delta^2$ would not cross
any other strand, a contradiction. Hence $a\neq 1$. In the same way
one shows that $a\neq n$. Therefore the permutation induced by $s$
can be written as
$$
 \pi_s= (a)(1 \ u_1 \ u_2 \ \cdots\  u_r \ n \ d_t \ d_{t-1} \ \cdots \ d_1).
$$
We can show that $u_1<\cdots<u_r$ and that $d_t>\cdots >d_1$, using
the same proof as in Proposition~\ref{P:USS of delta}. In
Figure~\ref{F:conj_and_notconj_of_epsilon} we can see an example of
two braids whose permutations are cycles of length $n-1$. The first
one satisfies the above inequalities and the second one does not.

\begin{figure}[ht]
\centerline{\includegraphics{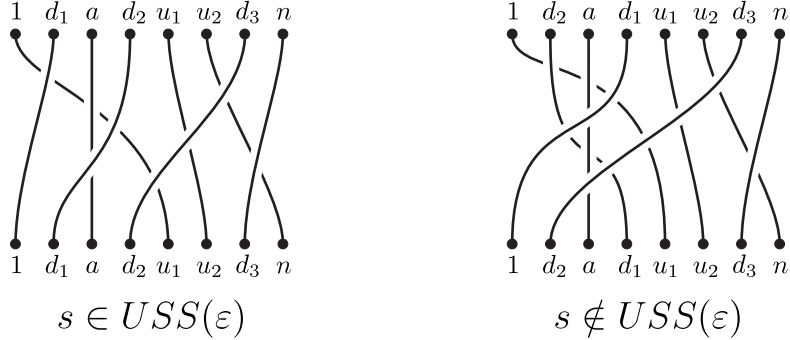}}
\caption{Two simple braids in $B_8^A$ whose permutations are cycles of
length $7$. By Proposition~\ref{P:USS of epsilon}, the first one is
conjugate to $\varepsilon$ and the second one is not. As in
Figure~\ref{F:conj_and_notconj_of_delta}, the exponent sums of the
two braids differ; the exponent sum of second one is 12, while the
exponent sum of conjugates of $\varepsilon\in B_8^A$ is 8.}
\label{F:conj_and_notconj_of_epsilon}
\end{figure}

Now let $s$ be a simple element such that
$$
 \pi_s= (a)(1 \ u_1 \ u_2 \ \cdots\  u_r \ n \ d_t \ d_{t-1} \ \cdots \ d_1)
$$
for some $u_1<u_2<\cdots < u_r$, some $d_t>d_{t-1}>\cdots >d_1$ and
some $a\neq 1$ or $n$.

Suppose that $t>0$. Similarly to the proof of Proposition~\ref{P:USS
of delta}, we will conjugate $s$ by $\sigma_{[d_1\rightarrow 1]}$,
and this will reduce the index $t$. Let $k=d_1$. If $a>k$ one has
$$
  \pi_s=(a)(\underline{1\ 2\ \cdots \ k-1} \ u_{k-1}\ \cdots \ u_r\ \ n \ d_t \ \cdots \
d_2 \underline{\ k}\:),
$$
otherwise
$$
  \pi_s=(a)(\underline{1\ 2\ \cdots \ a-1 \ a+1\ \cdots \ k-1} \ u_{k-2}\ \cdots \ u_r\ \ n \ d_t \ \cdots \
d_2 \underline{\ k}\:).
$$
The picture in the former case is the same as in
Figure~\ref{F:conjugate_of_delta}, while the latter case is
represented in Figure~\ref{F:conjugate_of_epsilon}. In either case,
the strand of $s$ that ends at $k$ (that is, $d_2$ if $t>1$ or $n$
if $t=1$) does not cross the strands that end at $1,2,\ldots,k-1$
(that is $k,1,2,\cdots,k-2$, where one of them could possibly be
equal to $a$). Therefore $s \:\sigma_{[k\rightarrow 1]}$ is simple.
At the same time, the strand $k$ of $s$ crosses the strands
$k-1,k-2,\ldots, 1$ (where one of them could be equal to $a$). Hence
$s'=\sigma_{[k\rightarrow 1]}^{-1} s \:\sigma_{[k\rightarrow 1]}$ is
simple. Depending on whether $a>k$ or not, one has either
$$
  \pi_{s'}=(a)(\underline{1\ 2\ \cdots \ k-1 \ k} \ u_{k-1}\ \cdots \ u_r\ \ n \ d_t \ \cdots \
d_2 ),
$$
or
$$
  \pi_{s'}=(a+1)(\underline{1\ 2\ \cdots \ a \ a+2\ \cdots \ k} \ u_{k-2}\ \cdots \ u_r\ \ n \ d_t \ \cdots \
d_2 ).
$$

\begin{figure}[ht]
\centerline{\includegraphics{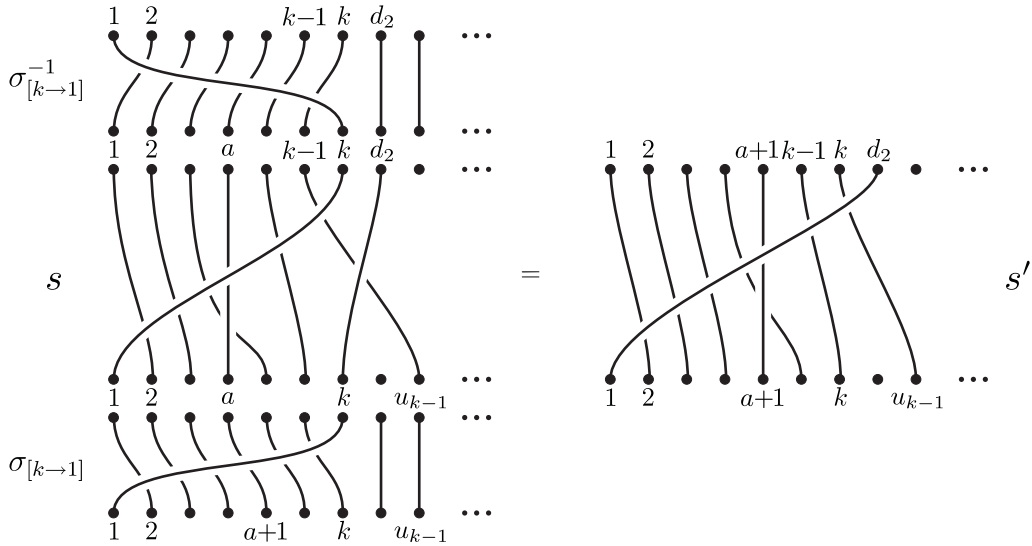}}
\caption{Conjugating $s$ to $s'$ when $a<k$.}
\label{F:conjugate_of_epsilon}
\end{figure}

We can continue this process, increasing the number of indices
between $1$ and $n$. Notice that the index $a$ increases at some
step $i$ if and only if $d_i>a$, and in that case we will also have
$d_{i+1},\ldots,d_t>a$, so once the index $a$ increases, it
continues increasing at every further step of the procedure.
Eventually, one obtains a simple element $s_0= \alpha^{-1} s \alpha
$ such that
$$
 \alpha= \sigma_{[d_1\rightarrow 1]}\sigma_{[d_2\rightarrow 1]}\cdots \sigma_{[d_t\rightarrow 1]}
$$
and
$$
\pi_{s_0} = (b)(1\ 2\ \cdots \ b-1 \ b+1 \ \cdots \ n),
$$
where $b=a+t-\max\{i\,:\,d_i<a \}$.

If $b=2$ we already have $s_0=\varepsilon$. Otherwise we will
conjugate $s_0$ by $\sigma_{[b\rightarrow 2]}$. Notice that the
strand of $s_0$ that ends at $b$ (that is, strand $b$ itself) does
not cross the strands ending at $b-1,b-2,\ldots,2$. Hence $s_0
\sigma_{[b\rightarrow 2]}$ is simple. Next, the strand $b$ of $s_0
\sigma_{[b\rightarrow 2]}$ crosses the strands $b-1,b-2,\ldots,2$,
hence $\sigma_{[b\rightarrow 2]}^{-1} s_0 \sigma_{[b\rightarrow 2]}$
is simple, and its permutation is equal to $(2)(1\ 3\ 4\ \cdots\
n)$, hence this simple braid is equal to $\varepsilon$. Therefore,
if we define
$$
 \beta = \sigma_{[d_1\rightarrow 1]}\:\sigma_{[d_2\rightarrow
1]}\cdots \sigma_{[d_t\rightarrow 1]}\: \sigma_{[b\rightarrow 2]}
$$
where $b=a+t-\max\{i\,:\,d_i<a \}$, then $\beta^{-1} s \beta =
\varepsilon$.
\end{proof}

\begin{remark} {\rm  The element $\beta$ defined above is not
necessarily simple, but in the worst case it is the product of two
simple elements, $\sigma_{[d_1\rightarrow
1]}\:\sigma_{[d_2\rightarrow 1]}\cdots \sigma_{[d_t\rightarrow 1]} $
and $\sigma_{[b\rightarrow 2]}$. Hence, every element in
$USS(\varepsilon)$ is connected to $\varepsilon $ by a conjugating
element of canonical length at most 2.}
\end{remark}

\begin{corollary}\label{C:size of USS(epsilon)}
If $\varepsilon=\sigma_1(\sigma_{n-1}\cdots \sigma_{1})\in B_n^A$ then
$\# (USS(\varepsilon))= (n-2)2^{n-3}$.
\end{corollary}

\begin{proof}
The elements in $USS(\varepsilon)$ are characterized by the
permutation given in Proposition~\ref{P:USS of epsilon}, which is
itself characterized by the sequence $1<u_1<\cdots <u_r<n$ and the
number $a$.  Since $a\neq 1,n$, one has $n-2$ choices for the index
$a$. And for every choice of $a$, the number of possible sequences
is equal to the number of subsets of $\{2,\ldots,a-1, a+1,
\ldots,n-1\}$, which is precisely $2^{n-3}$. Hence the total number
of choices is $(n-2)2^{n-3}$.
\end{proof}

Notice that the results in this section not only characterize the
elements in $USS(\delta)$ and $USS(\varepsilon)$ by their
permutations, determining the sizes of these two sets, but also find
conjugating elements from any given element in $USS(\delta)$ (resp.
$USS(\varepsilon)$) to $\delta$ (resp. $\varepsilon$). This fact,
together with the known algorithm for obtaining for any braid $X$ a
conjugate $Y$ of $X$ whose canonical length is minimal~\cite{El-M},
which for periodic $X$ implies $Y\in USS(X)$, provides a
solution to the conjugacy search problem for conjugates of $\delta$ or
$\varepsilon$. Moreover, this algorithm has complexity $O(l^3n^3\log
n)$, where $l$ is the letter length in Artin generators of the input
braid. But this algorithm is not easily generalized to other
periodic braids (conjugates of {\it powers} of $\delta$ or
$\varepsilon$), so in the next section we will present an
alternative approach that solves the conjugacy search problem for
every periodic braid, using other Garside structures and other
groups (namely Artin-Tits groups of type $\mathbf B$).

\begin{remark} {\rm  In~\cite{M-H} there is a simple algorithm
which finds a conjugating element from any braid in $USS(\delta)$ to
$\delta$. It is also an efficient algorithm, and very easy to
implement, but Proposition~\ref{P:USS of delta} directly provides a
conjugating element and at the same time characterizes the elements
in $USS(\delta)$, in such a way that we can count them all.}
\end{remark}

\begin{remark} {\rm
We end this section by remarking that, in practical computations for
small $n$, the sizes of $USS(\delta^k)$ and $USS(\varepsilon^k)$ for
different values of $k$ are in most cases much bigger than the sizes
of $USS(\delta)$ and $USS(\varepsilon)$, respectively. Hence the
usual algorithm in~\cite{Gebhardt} is not efficient in general for
periodic braids. We also notice that the algorithm in~\cite{M-H} can
be generalized to $\varepsilon$, but it does not generalize to
powers of $\delta$ or $\varepsilon$. Hence the algorithm in the next
section is, to our knowledge, the first efficient algorithm to solve
the conjugacy search problem for periodic braids.}
\end{remark}

\section{Proof of Theorem~\ref{T:main theorem}}
\label{S:proof}

In this section we will complete the proof of Theorem~\ref{T:main
theorem} by developing a polynomial algorithm to solve the conjugacy
search problem for periodic braids.

Suppose that we are given two braids $X,Y\in B_n^A$. Using Algorithm
A, we may assume that $X$ and $Y$ are periodic, and that they are
conjugate to the same power of $\delta$ or $\varepsilon$ (otherwise
we would stop and return a negative answer for steps (1) or (2) in
Theorem~\ref{T:main theorem}). We can also assume that we know the
specific power of $\delta$ (resp. $\varepsilon$) which is conjugate
to $X$ and $Y$, say $\delta^k$ (resp. $\varepsilon^k$). Clearly, we
just need an algorithm that finds a conjugating element from $X$ to
$\delta^k$ (resp. $\varepsilon^k$), since the same algorithm can be
applied to $Y$ and we would immediately obtain a conjugating element
from $X$ to $Y$.

Therefore, we will suppose that we are given a braid $X\in B_n^A$ as a
word of length $l$ in $\sigma_1,\ldots, \sigma_{n-1}$ and their
inverses, and that $X$ is conjugate to $\delta^k$ respectively
$\varepsilon^k$ for some $k\neq 0$. We will describe algorithms
finding a conjugating element from $X$ to $\delta^k$ or
$\varepsilon^k$, whose complexities are polynomial in $n$ and
$l$. The two cases are treated separately, in Sections \ref{SS:powers of delta} and \ref{SS:powers of varepsilon} below.

\subsection{Solving the conjugacy search problem for conjugates of $\delta^k$}
\label{SS:powers of delta}

\begin{figure}[ht]
\centerline{\includegraphics{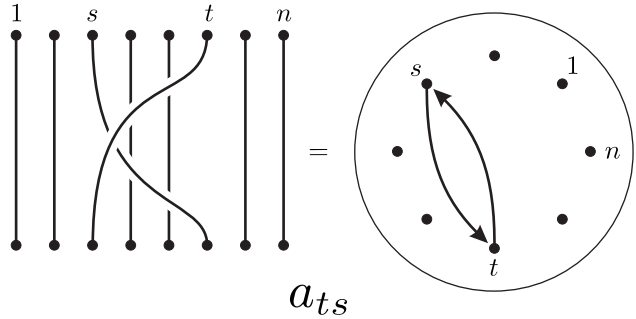}} \caption{The band
generator $a_{ts}$.} \label{F:band_generators}
\end{figure}

We remind the reader that in~\cite{BKL}, Birman, Ko and Lee
investigated a then-new presentation for the braid groups:
\begin{equation}
\label{E:BKL presentation} B_n^{B}:  \left< a_{ts}, \ \ \mbox{\scriptsize
$1\leq s < t \leq n$} \left|
\begin{array}{ll}
a_{ts}a_{rq} = a_{rq}a_{ts}   & \mbox{\scriptsize if \; $(t-r)(t-q)(s-r)(s-q)>0$},  \\
a_{ts}a_{sr} = a_{tr}a_{ts} = a_{sr}a_{tr} &  \mbox{\scriptsize if
\; $1\leq r< s < t \leq n$.}\ \end{array}\right.
 \right>.
\end{equation}
The elements $a_{ts}$ are called {\it band generators} or {\it
Birman-Ko-Lee generators}. The left sketch in Figure ~\ref{F:band_generators} shows one way to think of the generator $a_{ts}$.  A different way is shown on the right, where we consider $\mathbb D^2_n$ to be the disc in
$\mathbb C$ centered at the origin with radius 2, the $n$ punctures
being the $n$-th roots of unity $\zeta_k=e^{2k\pi i/n}$ for
$k=1,\ldots, n$.  Then $a_{ts}$ is the braid that swaps the punctures
$\zeta_s$ and $\zeta_t$ as shown in the right hand side of
Figure~\ref{F:band_generators}.

In will be convenient to think of $B_n^{A}$ and $B_n^{B}$ as defining distinct groups. The relation between them is then given by the isomorphism $\Phi:B_n^{B} \to B_n^{A} $:
\begin{eqnarray}\label{E:the isomorphism BKL to Artin}
\Phi(a_{i,i+1}) & = &\sigma_i, \ \ 1\leq i\leq n-1\cr
\Phi(a_{ts}) & =&  (\sigma_{t-1} \sigma_{t-2} \cdots \sigma_{s+1}) \sigma_s (\sigma_{s+1}^{-1} \cdots
\sigma_{t-2}^{-1}\sigma_{t-1}^{-1}), \ \ 1\leq s< t-1 \leq n.
\end{eqnarray}
The inverse automorphism sends $\sigma_i$ to $a_{i,i+1}$.

The reason we wish to think of these two groups as being distinct, is because we need to distinguish the Garside structure on $B_n^A$ \cite{Garside} from that on $B_n^B$, introduced in \cite{BKL}.  When we say that $X\in B_n^{A}$ (resp.
$X\in B_n^{B}$) is written in left normal form, our notation is intended to mean that we are using the Garside structure associated to the presentation (\ref{E:classical presentation}) (resp. (\ref{E:BKL presentation})).  The key point here (which we will generalize when we treat the case of braids conjugate to $\varepsilon$, is that the Garside element for $B_n^B$ is
precisely our periodic braid  $\delta$.   It is shown in~\cite{BKL} that with respect to the Garside structure introduced in \cite{BKL}, the left normal form of a braid in $B_n^B$, given as a word
of length $l$ in the band generators and their inverses, can be
computed in time $O(l^2n)$.  We will solve the conjugacy search problem for braid conjugate to $\delta$ by making use of the algorithm in \cite{Gebhardt}, using the Garside structure on $B_n^B$.  This will enable us to bypass the difficulty which was uncovered in Corollary~\ref{C:size of USS(epsilon)}. 

\begin{figure}[ht]
\centerline{\includegraphics{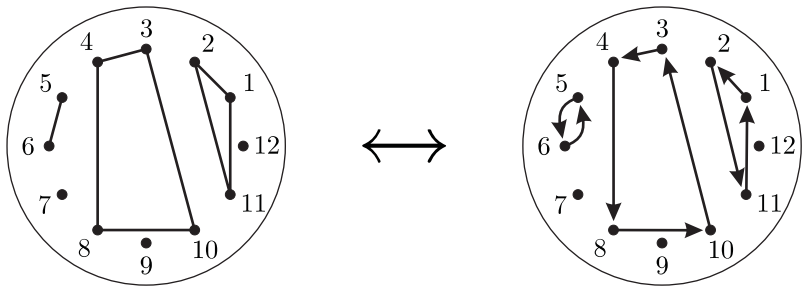}} \caption{A simple
element of $B^B_{12}$, which is a
product of three polygonal braids. It is the braid
$(a_{11,2}a_{2,1})(a_{10,8}a_{8,4}a_{4,3})(a_{6,5})$, where the
three factors (the polygonal braids) commute.}
\label{F:simple_elements}
\end{figure}

It will be important for our purposes to describe the simple
elements in the Garside structure on $B_n^B$.  These simple
elements are known to be in bijection with the {\it non-crossing
partitions} of the $n$-th roots of unity $\mathcal
R=\{\zeta_1,\ldots,\zeta_n\}$~\cite{BKL,Bessis}. Non-crossing
partitions can be defined as follows: Given a partition $\wp$ of
$\mathcal R$, every part of $\wp$ with $d$ elements ($d\geq 2$)
gives rise to a unique convex polygon joining the $d$ punctures (if
$d=2$ the polygon is just a segment). The partition $\wp$ is said to
be non-crossing if these polygons are pairwise disjoint.  Each
polygon determines a braid which corresponds to a rotation of its
$d$ vertices in the counterclockwise sense, and that we will call a
{\it polygonal braid}. Disjoint polygons determine commuting
polygonal braids. The simple element corresponding to a non-crossing
partition $\wp$ is the product of the (mutually commuting) polygonal
braids determined by $\wp$, as is shown in
Figure~\ref{F:simple_elements}. Hence, each simple element of
$B_{n}^B$ is a product of at most $n/2$ polygonal braids. Notice also
that the polygonal braid corresponding to the part
$\{\zeta_{i_1},\zeta_{i_2}, \cdots ,\zeta_{i_k}\}$, with $i_1<i_2<
\cdots < i_k$, is precisely $a_{i_k,i_{k-1}} \: a_{i_{k-1},
i_{k-2}}\cdots a_{i_2,i_1}$. The element $\delta$ is the polygonal
braid corresponding to the whole set $\{\zeta_1,\ldots,\zeta_n\}$.

Before stating and proving the main result of this section, we need
a lemma that will improve the estimation of the complexity of our
algorithms for periodic braids. It is the following:

\begin{lemma}\label{L:n<l}
If a nontrivial periodic braid $X\in B_n^A$ is given as a word of
length $l$ in the Artin generators and their inverses, then $l\geq
n-1$.
\end{lemma}

\begin{proof}
Suppose that $l<n-1$. Then the exponent sum of $X$ is an integer $m$
with $0\leq |m| < n-1$. By Lemma~\ref{L:exponent sum}, the exponent
sum of a periodic braid is a multiple of either $n-1$ or $n$. It
follows that $m=0$, so $X$ is conjugate to $\delta^0=1$. But in this
case $X$ is trivial, a contradiction.
\end{proof}

We can finally show our main result for conjugates of powers of
$\delta$.

\begin{proposition}\label{P:Algorithm B}
Let $X\in B_n^A$ be given as a word of length $l$ in the Artin
generators $\sigma_1,\ldots,\sigma_{n-1}$ and their inverses. If $X$
is conjugate to $\delta^k$ for $k\neq 0$, there exists an algorithm
of complexity $O(l^3n^2)$ that finds a conjugating element $C\in
B_n^A$ such that $C^{-1}XC=\delta^k$.
\end{proposition}

\begin{proof}
We are given a word $X\in B_n^A$:
$$
 X = \sigma_{\mu_1}^{\epsilon_1}\sigma_{\mu_2}^{\epsilon_2}\cdots\sigma_{\mu_l}^{\epsilon_l}.
$$

It is very simple to rewrite $X$ as a word in the band generators,
because $\Phi^{-1}(\sigma_i)=a_{i+1,i}$ for each $i=1\dots,n-1$.  So we have:
$$
 \Phi^{-1}(X) = a_{\mu_1+1,\mu_1}^{\epsilon_1}a_{\mu_2+1,\mu_2}^{\epsilon_2}\cdots
a_{\mu_l+1,\mu_l}^{\epsilon_l}.
$$
We can then apply iterated cycling and decycling to $X\in
B_n^B$, in order to obtain a conjugate $X'\in B_n^B$ of
minimal canonical length, together with a conjugating element. It is
shown in~\cite{BKL-2001} that we need to apply at most $|\delta|\,
l$ cyclings and decyclings, this means at most $|\delta|\, l$
computations of normal forms, where $|\delta|$ is the letter length
of $\delta$ written as a positive word in the band generators. Since
$|\delta|=n-1$, it follows that we can obtain $X'\in B_n^B$ of
minimal canonical length, and a conjugating element from $X$ to
$X'$, in time $O(l^3n^2)$.

But $X$ is conjugate to $\delta^k$, which is a power of the Garside
element of $B_n^B$, so $\ell(\delta^k)=0$. This means that
$USS(X)=\{\delta^k\}$, and more precisely $X'=\delta^k$. Hence, we
have found a conjugating element $C\in B_n^B$ from $X$ to
$\delta^k$ in time $O(l^3n^2)$. As this conjugating element is given
in terms of band generators, the last step consists of translating
$C$ to Artin generators. 

Recall that a cycling (resp. a decycling) consists of a conjugation
by a simple element (resp. by the inverse of a simple element). So
$C$ is a product of at most $(n-1)l$ simple elements (or inverses)
in $B_n^B$. In the Birman-Ko-Lee structure, the letter length
of a simple element is at most $n-1$, so $C\in B_n^B$ has
letter length at most $(n-1)^2 l$. Since each band generator is
equal to a word in Artin generators of length at most $2n-3$, this
means that one can translate $C$ to Artin generators, via the isomorphism $\Phi$, in time
$O(n^3l)$.

Therefore, the conjugacy search problem for conjugates of
$\delta^k$, given as words in Artin generators, can be solved in
time $O(l^3n^2+ln^3)$.  By Lemma~\ref{L:n<l}, one has $l\geq
n-1$, so that $ln^3\leq l(l+1)n^2 < l^3n^2$ (we can assume $l>1$).
Hence this complexity is equal to $O(l^3n^2)$.
\end{proof}

The algorithm described in the proof of
Proposition~\ref{P:Algorithm B} is the following.

\medskip
\noindent {\bf Algorithm B:}

Input: A word $w$ in Artin generators and their inverses
representing $X\in B_n^A$ conjugate to $\delta^k$.

Output: $C\in B_n^A$ such that $C^{-1}XC=\delta^k$.

\begin{enumerate}

 \item Translate $w$ to a word $w'$ in band generators using the rule $\sigma_i
\rightarrow a_{i+1,i}$.

 \item Apply iterated cyclings and decyclings in $B_n^B$ to $w'$
   until $\delta^k$ is obtained. Let $C'\in B_n^B$ be the product of
all the conjugating elements in this process.

 \item Translate $C'$ to a word $C\in B_n^A,$ using the
rule
$$
a_{ts} \rightarrow (\sigma_{t-1} \sigma_{t-2} \cdots \sigma_{s+1})
\sigma_s (\sigma_{s+1}^{-1} \cdots
\sigma_{t-2}^{-1}\sigma_{t-1}^{-1}).
$$

\item Return $C$.

\end{enumerate}

\medskip

By Proposition~\ref{P:Algorithm B}, Algorithm $B$ has complexity
$O(l^3n^2)$.

\subsection{Solving the conjugacy search problem for conjugates of $\varepsilon^k$}
\label{SS:powers of varepsilon}

Our final task is to learn how to find the conjugating element in
the case when $X$ is conjugate to $\varepsilon^k$.   The methods will be identical to those used in case of conjugates of $\delta^k$:  We will begin with $X\in B_n^A$, i.e. $X$ will be given as a word in the generators of $B_n^A$ and their inverses.  Using Algorithm A we will have verified that $X$ is periodic and conjugate to a known power of $\varepsilon$.  Our task will be to find the conjugating element.  
We
will prove that there is also a suitable Garside group, with a known Garside structure, whose
Garside element is $\varepsilon$.  This group, however,  is not the braid group, rather it is a subgroup of the braid group that we will denote $P_{n,2}$. The subgroup is formed by the
braids whose corresponding permutation preserves the second
puncture. It is well known that $P_{n,2}$ is a Garside group, since
it is isomorphic to the Artin-Tits group of type $\mathbf
B_{n-1}$~\cite{Crisp}.  Nevertheless, we won't use the classical
Garside structure on the Artin-Tits group $\mathcal A(\mathbf B_{n-1})$, but the {\it
dual} Garside structure defined in~\cite{Bessis}.  This explains why we
shall start, in Section~\ref{SS-3.1.1}, by describing the groups, embeddings and Garside
structures that we will need to use in our algorithm.   We then put them to work in Section~\ref{SS-3.1.2}

\subsubsection{Braids fixing one puncture, Artin-Tits groups of type
  $\mathbf B$ and symmetric braids} \label{SS-3.1.1}

We shall now describe the five groups we are interested in, with
their corresponding Garside structures. The first two groups are
well known, they are just $B_n^A$ and $B_{2n-2}^B$.

Next, let us consider the subgroup $P_{n,2}\subset B_n^A$,
consisting of braids that fix the second puncture. That is,
$P_{n,2}=\{X\in B_n^A\,:\,\pi_X(2)=2\}$. We will not consider right now
a Garside structure on $P_{n,2}$, but we remark that $\varepsilon\in
P_{n,2}$.

Now let $Sym_{2n-2}$ be the centralizer of $\delta^{n-1}$ in
$B_{2n-2}^B$, where we write $\delta$ for $\delta_{2n-2}$. In other words, if we represent the $2n-2$
punctures of $\mathbb D^2_{2n-2}$ as the $(2n-2)$-nd roots of unity,
the elements of $Sym_{2n-2}$ are precisely the braids which are
invariant under a rotation of 180 degrees. This is why they are
called {\it symmetric braids}.

Finally, consider the Artin-Tits group $\mathcal A(\mathbf
B_{n-1})$, whose presentation is
$$
 \mathcal A(\mathbf B_{n-1})=\left\langle s_1,\ldots, s_{n-1} \; \left|
\begin{array}{ll}
 s_is_j= s_js_i & \mbox{ if } |i-j|>1 \\
 s_is_js_i=s_js_is_j & \mbox{ if } |i-j|=1 \mbox{ and } i,j\neq 1
\\
s_1s_2s_1s_2=s_2s_1s_2s_1
\end{array}
  \right. \right\rangle.
$$

We shall now recall from the literature that the last three groups
we just considered are isomorphic, that is, one has the following
situation:
\begin{eqnarray} \label{E:the two Garside structures}
  \begin{array}{ccccc}
    B_n^A & & & & B_{2n-2}^B \\
        \cup      & & & &   \cup   \\
   P_{n,2}   & \cong & \mathcal A(\mathbf B_{n-1}) & \cong &
Sym_{2n-2}
 \end{array}
\end{eqnarray}
Moreover, it can be deduced from~\cite{Reiner} that the restriction
of the Garside structure of $B_{2n-2}^B$ determines a Garside
structure in $Sym_{2n-2}$. Via the above isomorphisms, this
induces Garside structures in $\mathcal A(\mathbf B_{n-1})$ and
in $P_{n,2}$. We shall see that the Garside element of the latter is
precisely $\varepsilon$, and this will help us to solve the
conjugacy search problem for conjugates of $\varepsilon^k$.

Let us then study in detail the mentioned isomorphisms.

\begin{lemma}\label{L:from Artin B to P}
The map $\rho:\:\mathcal A(\mathbf B_{n-1}) \rightarrow P_{n,2}$
given by $\rho(s_1)=\sigma_1^2$,
$\rho(s_2)=\sigma_1\sigma_2\sigma_1^{-1}$ and $\rho(s_i)=\sigma_i$
for $i>2$, is an isomorphism.
\end{lemma}

\begin{proof}
Proposition 5.1 in~\cite{Crisp} provides an isomorphism $\rho_0:\
\mathcal A(\mathbf B_{n-1})\rightarrow P_{n,1}$, where $P_{n,1}$ is
the subgroup of $B_n^A$ consisting of braids which fix the first
puncture. This isomorphism is given by $\rho_0(s_1)=\sigma_1^2$ and
$\rho_0(s_i)=\sigma_i$ for $i>1$, and it was already known to
specialists, prior to~\cite{Crisp}. Now we just need to notice that
the inner automorphism $\varphi: B_n^A \rightarrow B_n^A$ given by
$\varphi(X)=\sigma_1 X\sigma_1^{-1}$ sends $P_{n,1}$ isomorphically
to $P_{n,2}$, and that $\varphi_{|_{P_{n,1}}}\circ \rho_0 = \rho$.
\end{proof}

\begin{figure}[ht]
\centerline{\includegraphics{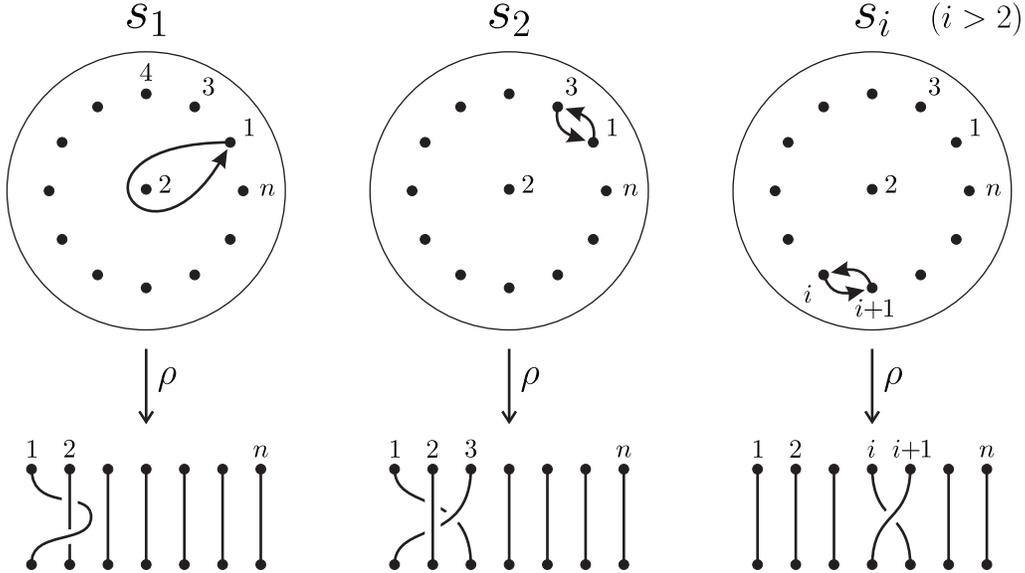}}
\caption{The generators of $\mathcal A(\mathbf B_{n-1})$,
represented as braids on $\mathbb D^2\backslash\{0\}$, and their
images under the isomorphism $\rho:\ \mathcal A(\mathbf B_{n-1})
\rightarrow P_{n,2}$.} \label{F:Artin_Tits_generators}
\end{figure}

\begin{remark}{\rm It is well known~\cite{Crisp} that $P_{n,2}$ (hence $\mathcal
A(\mathbf B_{n-1})$) can be identified with the braid group of the
open annulus $\mathbb D^2\backslash \{0\}$ on $n-1$ strands. Indeed,
an element $X\in P_{n,2}$ fixes the second puncture, so it can be
isotoped to a braid whose second strand in $\mathbb D^2\times [0,1]$
is a straight line, say $\{0\}\times [0,1]$. This second strand can
be considered to be a hole of $\mathbb D^2$, so $X$ can be regarded
as a braid on $n-1$ strands of $\mathbb D^2\backslash \{0\}$.

In order to avoid confusion, we will represent elements in
$P_{n,2}\in B_n^A$ in the usual way, as they are represented
at the bottom of Figure~\ref{F:Artin_Tits_generators}, while
elements of $\mathcal A(\mathbf B_{n-1})$ will be represented in the
Birman-Ko-Lee style, as braids on $\mathbb D^2\backslash\{0\}$ whose
base points are the $(n-1)$-st roots of unity, as we can see at the
top of Figure~\ref{F:Artin_Tits_generators}.  }
\end{remark}

\begin{lemma}\label{L:from Artin B to Sym}
The map $\theta':\:\mathcal A(\mathbf B_{n-1}) \rightarrow
Sym_{2n-2}$ given by $\theta'(s_1)=a_{n,1}$ and
$\theta'(s_i)=a_{i,i-1}\, a_{i+n-1,i+n-2}$ for $i>1$, is an
isomorphism.
\end{lemma}

\begin{figure}[ht]
\centerline{\includegraphics{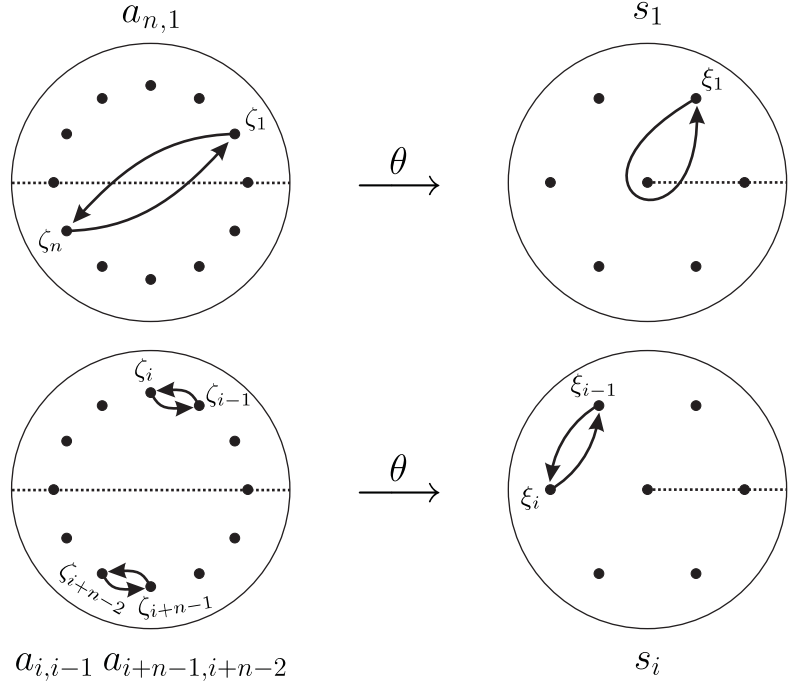}} \caption{The map $\theta$
transforms the symmetric braids on the left hand side to the
generators of $\mathcal A(\mathbf B_{n-1})$ on the right hand side.}
\label{F:theta}
\end{figure}

\begin{proof}
In~\cite{Brieskorn}, Brieskorn showed that an Artin-Tits group of
finite type is the fundamental group of the regular orbit space of
its corresponding Coxeter group, acting as a finite real reflection
group on a complex space. In particular, since the Coxeter group
associated to $\mathcal A(\mathbf B_{n-1})$ is
$W=\Sigma_{n-1}\ltimes (\mathbb Z/ 2 \mathbb Z )^{n-1}$, where the
symmetric group acts by permuting coordinates (that is, $W$ is the
signed permutation group), and its corresponding hyperplane
arrangement is $x_1 x_2 \cdots x_{n-1}\prod_{i\neq
j}(x_i-x_j)(x_i+x_j)$, it follows that $\mathcal A(\mathbf
B_{n-1})=\pi_1(X_{\mathbf B_{n-1}}/W)$, where
$$
X_{\mathbf B_{n-1}}=\{(x_1,\ldots,x_{n-1})\in \mathbb C^{n-1} \ | \
x_i\neq \pm x_j \mbox{ for } i\neq j; \; x_i\neq 0 \mbox{ for all }
i \}.
$$
A good way to describe the space $X_{\mathbf B_{n-1}}$ is as the set
of $(n-1)$-tuples of pairs
$$
((x_1,-x_1),(x_2,-x_2),\ldots, (x_{n-1},-x_{n-1})),
$$
where each $x_i\in\mathbb C$, any two pairs are distinct, and
$x_i\neq 0$ for all $i$. Considering the action of $W$, all the above
pairs and $(n-1)$-tuples can be regarded as unordered. Hence
$X_{\mathbf B_{n-1}}/W$ is
the configuration space of $2n-2$ disjoint and undistinguishable
points in $\mathbb C$, whose configuration is invariant under
multiplication by $-1$.  We can choose as a base point of this space
the $(2n-2)$-nd roots of unity. Hence, an element of its fundamental
group is represented by a braid which is invariant under a rotation
by 180 degrees, that is, by a symmetric braid in $B_{2n-2}^B$.

It is important to note that two symmetric braids represent the same
element in $\pi_1(X_{\mathbf B_{n-1}}/W)$ if and only if they are
isotopic {\it through symmetric braids}, hence one cannot say a
priori that two symmetric braids that are isotopic in
$B_{2n-2}^B$ represent the same element of $\pi_1(X_{\mathbf
B_{n-1}}/W)$. Fortunately, it is shown in~\cite{BDM} that two
symmetric braids are isotopic in $B_{2n-2}^B$ if and only if
they are isotopic through symmetric braids. That is, it is shown
that $\mathcal A(\mathbf B_{n-1})= \pi_1(X_{\mathbf B_{n-1}}/W)
\cong Sym_{2n-2}$.

Moreover, from the work in~\cite{BDM} one obtains an isomorphism
$\theta: \ Sym_{2n-2} \rightarrow \mathcal A(\mathbf B_{n-1})$,
where elements of $Sym_{2n-2}$ are symmetric braids based on the
$(2n-2)$-nd roots of unity, and the elements of $\mathcal A(\mathbf
B_{n-1})$ are considered as braids on the annulus $\mathbb
D^2\backslash\{0\}$ based on the $(n-1)$-st roots of unity. The
isomorphism $\theta$ can be easily described geometrically, since it
just identifies antipodal points in $\mathbb C$. That is, it sends
$z\in \mathbb C\backslash \{0\}$ to $z^2/|z|$.  This corresponds to
a two-sheeted covering map of $\mathbb C\backslash \{0\}$, and since
no strand of a symmetric braid touches the axis $\{0\}\times [0,1]$,
this map is well defined.

In Figure~\ref{F:theta} we can see that $\theta(a_{n,1})=s_1$ and
that $\theta(a_{i,i-1}\, a_{i+n-1,i+n-2})=s_i$ for $i>1$, where in
the picture one has $\zeta_k=e^{2k \pi i/(2n-2)}$ and
$\xi_k=e^{2k\pi i/(n-1)}$. Therefore $\theta'=\theta^{-1}$, so it is
an isomorphism.
\end{proof}

By Lemmas~\ref{L:from Artin B to P} and \ref{L:from Artin B to Sym}
we know that $P_{n,2} \cong \mathcal A(\mathbf B_{n-1}) \cong
Sym_{2n-2}$, and we also know how to transform any word in the
generators $s_1,\ldots,s_{n-1}$ of $\mathcal A(\mathbf B_{n-1})$ and
their inverses, into a word in either the Artin generators of
$P_{n,2}$ or the band generators of $Sym_{2n-2}$, via the
isomorphisms $\rho$ and $\theta'=\theta^{-1}$.
$$
  \begin{array}{ccccc}
    B_n^A & & & & B_{2n-2}^B \\
        \cup      & & & &   \cup   \\
   P_{n,2}   & \stackrel{\rho}{\longleftarrow} & \mathcal A(\mathbf B_{n-1}) & \stackrel{\theta'}{\longrightarrow} &
Sym_{2n-2}.
 \end{array}
$$
 But in our algorithm we will
need to translate any word in the Artin generators of $B_n^A$,
representing an element of $P_{n,2}$, to a word in the band
generators of $Sym_{2n-2}$, and vice versa. Hence, we need the
following results.

\begin{lemma}\label{L:from P to Sym}
Let $X\in P_{n,2}\subset B_n^A$ be given as a word of length
$l$ in the Artin generators and their inverses,
$X=\sigma_{\mu_1}^{\epsilon_1}\sigma_{\mu_2}^{\epsilon_2}\cdots
\sigma_{\mu_l}^{\epsilon_l}$. For $i=0,\ldots,l,$ let
$X_i=\sigma_{\mu_1}^{\epsilon_1}\sigma_{\mu_2}^{\epsilon_2}\cdots
\sigma_{\mu_i}^{\epsilon_i}$ and let $k_i=\pi_{X_i}(2)$, that is,
the final position of the second strand of $X_i$. Then one obtains a
word in the band generators and their inverses representing
$\theta'(\rho^{-1}(X))\in Sym_{2n-2}$, by replacing each letter
$\sigma_{\mu_i}^{\epsilon_i}$ using the following rules:
$$
\sigma_{\mu_i} \rightarrow \left\{
\begin{array}{ll}
  a_{\mu_i+1,\mu_i}\: a_{\mu_i+n,\mu_i+n-1} & \mbox{ if} \quad \mu_i < k_{i-1}- 1, \medskip \\
  1  & \mbox{ if} \quad\mu_i=k_{i-1}-1, \medskip \\
  a_{\mu_i+n-1,\mu_i} & \mbox{ if}\quad \mu_i=k_{i-1}, \medskip \\
  a_{\mu_i,\mu_i-1}\: a_{\mu_i+n-1,\mu_i+n-2} & \mbox{ if}\quad
\mu_i>k_{i-1},
\end{array} \right.
$$
and
$$
\sigma_{\mu_i}^{-1} \rightarrow \left\{
\begin{array}{ll}
  a_{\mu_i+n,\mu_i+n-1}^{-1}\: a_{\mu_i+1,\mu_i}^{-1}  & \mbox{ if} \quad \mu_i < k_{i-1}-
1, \medskip
\\
  a_{\mu_i+n-1,\mu_i}^{-1} & \mbox{ if} \quad\mu_i=k_{i-1}-1, \medskip \\
  1 & \mbox{ if}\quad \mu_i=k_{i-1}, \medskip \\
  a_{\mu_i+n-1,\mu_i+n-2}^{-1}\: a_{\mu_i,\mu_i-1}^{-1} & \mbox{ if}\quad
\mu_i>k_{i-1}.
\end{array} \right.
$$
Moreover, this algorithm has complexity $O(l)$, and produces a word
of length at most $2l$.
\end{lemma}

\begin{proof}
Recall that we are given a braid $X\in B_n^A$ that fixes the second
puncture, that is, $X\in P_{n,2}$, written as a word in the Artin
generators of $B_n^A$ and their inverses. We want to write
$\rho^{-1}(X)$ as a word in the generators $s_1,\ldots,s_{n-1}$ and
their inverses, and then $\theta'(\rho^{-1}(X))$ as a word in the
band generators of $B_{2n-2}^B$.

The first problem is that $X$ is not given as a word in the
generators of $P_{n,2}$, but in the generators of $B_n^A$.
We will then use the Reidemeister-Schreier method (see Section 2.3
of~\cite{MKS}) to decompose $X$ as a product of elements in
$P_{n,2}$. In order to do this, notice that $P_{n,2}$ is a subgroup
of $B_n^A$ of index $n$. The right coset of a braid $Z$ depends on where it
sends the second puncture. If $\pi_Z(2)=k$, we denote by $R_k$  a
representative of the right coset $P_{n,2}\,Z \in P_{n,2} \backslash B_n^A$. For
technical reasons, we will choose as coset representatives the
elements $R_1=\sigma_1$, $R_2=1$ and $R_k=\sigma_{[k\rightarrow
2]}^{-1}=\sigma_2^{-1}\cdots \sigma_{k-1}^{-1}$ if $k>2$.

Then, for $i=0,\ldots,l,$ we define $\overline{X_i}=R_{k_i}$. That
is, $\overline{X_i}$ is the chosen representative of $P_{n,2}\,X_i \in
P_{n,2} \backslash B_n^A$. Note that $\overline{X_0} = \overline{X_l}
= R_2 = 1$.

By the Reidemeister-Schreier method, one has
$$
  X=\prod_{i=1}^{l}\left(\overline{X_{i-1}}\: \sigma_{\mu_i}^{\epsilon_i}
\: \overline{X_i}\:^{-1}\right)= \prod_{i=1}^l \left(R_{k_{i-1}}
\sigma_{\mu_i}^{\epsilon_i} R_{k_i}^{-1}\right),
$$
where each of the above $l$ factors belongs to $P_{n,2}$. Notice
that $k_i=k_{i-1}$, unless either $\mu_i=k_{i-1}$ (in which case
$k_i=k_{i-1}+1$) or $\mu_i=k_{i-1}-1$ (and then $k_i=k_{i-1}-1$).
One can check that, depending on $\mu_i$ and $k_{i-1}$, each
of the above factors can be written in terms of the Artin generators
and their inverses as follows. If $\epsilon_i=1$, one has:
$$
(R_{k_{i-1}} \sigma_{\mu_i} R_{k_i}^{-1}) = \left\{
\begin{array}{ll}
  \sigma_2^{-1}\sigma_1 \sigma_2 & \mbox{ if} \quad 1=\mu_i < k_{i-1}- 1, \medskip \\
  \sigma_{\mu_i+1} & \mbox{ if} \quad 1\neq \mu_i < k_{i-1}- 1, \medskip \\
  1  & \mbox{ if } \hspace{1cm} \mu_i=k_{i-1}-1, \medskip \\
 \sigma_1^2 & \mbox{ if}\quad 1=\mu_i=k_{i-1}, \medskip \\
(\sigma_2^{-1}\sigma_3^{-1} \cdots \sigma_{\mu_i-1}^{-1})
\sigma_{\mu_i} (\sigma_{\mu_i} \sigma_{\mu_i-1}\cdots \sigma_2)  &
\mbox{
if}\quad 1\neq \mu_i=k_{i-1}, \medskip \\
 \sigma_1 \sigma_2 \sigma_1^{-1}  & \mbox{ if} \quad 2= \mu_i>k_{i-1}, \medskip \\
 \sigma_{\mu_i} & \mbox{ if}\quad 2\neq \mu_i>k_{i-1}.
\end{array} \right.
$$
If $\epsilon_i=-1$, one obtains the inverses of the above, in the
following way:
$$
(R_{k_{i-1}} \sigma_{\mu_i}^{-1} R_{k_i}^{-1}) = \left\{
\begin{array}{ll}
  \sigma_2^{-1}\sigma_1^{-1} \sigma_2 & \mbox{ if} \quad 1=\mu_i < k_{i-1}- 1, \medskip \\
  \sigma_{\mu_i+1}^{-1} & \mbox{ if} \quad 1\neq \mu_i < k_{i-1}- 1, \medskip \\
  \sigma_1^{-2} & \mbox{ if}\quad 1=\mu_i=k_{i-1}-1, \medskip \\
(\sigma_2^{-1}\sigma_3^{-1} \cdots \sigma_{\mu_i}^{-1})
\sigma_{\mu_i}^{-1}
(\sigma_{\mu_i-1}\cdots \sigma_2)  & \mbox{ if} \quad  1\neq \mu_i=k_{i-1}-1, \medskip \\
  1 & \mbox{ if }\hspace{1cm} \mu_i=k_{i-1}, \medskip \\
 \sigma_1\sigma_2^{-1} \sigma_1^{-1} & \mbox{ if}\quad 2=\mu_i>k_{i-1}, \medskip \\
 \sigma_{\mu_i}^{-1} & \mbox{ if}\quad 2\neq \mu_i>k_{i-1}.
\end{array} \right.
$$

Now we need to apply $\rho^{-1}$ to each factor $(R_{k_{i-1}}
\sigma_{\mu_i}^{\epsilon_i} R_{k_i}^{-1})$, and write the image in
in terms of the generators $s_1,\ldots, s_{n-1}$ and their inverses.
Recall that $\rho(s_1)=\sigma_1^2$, $\rho(s_2)=\sigma_1 \sigma_2
\sigma_1^{-1}= \sigma_2^{-1} \sigma_1 \sigma_2$ and
$\rho(s_i)=\sigma_i$ for $i>2$. Notice also that $\rho(s_2 s_1
s_2^{-1})=\sigma_2^2$, and that if $\mu_i>2$ one has
$$
  \rho\left((s_{\mu_i}s_{\mu_i-1} \cdots s_3\, s_2)\, s_1\, (s_2^{-1} s_3^{-1} \cdots
s_{\mu_i}^{-1})\right)= (\sigma_{\mu_i} \sigma_{\mu_i-1} \cdots
\sigma_3) \sigma_2^2 (\sigma_3^{-1} \cdots \sigma_{\mu_i}^{-1})
$$
$$
   = (\sigma_2^{-1} \cdots \sigma_{\mu_i-1}^{-1}) \sigma_{\mu_i}
(\sigma_{\mu_i}\cdots \sigma_2).
$$
Therefore, if $\epsilon_i=1$, one has:
$$
\rho^{-1}(R_{k_{i-1}} \sigma_{\mu_i} R_{k_i}^{-1}) = \left\{
\begin{array}{ll}
  s_{\mu_i+1} & \mbox{ if} \quad \mu_i < k_{i-1}- 1, \\
  1  & \mbox{ if} \quad\mu_i=k_{i-1}-1, \\
  (s_{\mu_i}s_{\mu_i-1} \cdots s_2)\: s_1\: (s_2^{-1} s_3^{-1}\cdots  s_{\mu_i}^{-1}) & \mbox{
if}\quad \mu_i=k_{i-1}, \\
 s_{\mu_i} & \mbox{ if}\quad \mu_i>k_{i-1},
\end{array} \right.
$$
and if $\epsilon_i=-1$, one obtains:
$$
\rho^{-1}(R_{k_{i-1}} \sigma_{\mu_i}^{-1} R_{k_i}^{-1}) = \left\{
\begin{array}{ll}
  s_{\mu_i+1}^{-1} & \mbox{ if} \quad \mu_i < k_{i-1}- 1, \\
   (s_{\mu_i}s_{\mu_i-1} \cdots s_2)\: s_1^{-1}\: (s_2^{-1}s_3^{-1}\cdots  s_{\mu_i}^{-1})  & \mbox{ if} \quad\mu_i=k_{i-1}-1, \\
  1 & \mbox{ if}\quad \mu_i=k_{i-1}, \\
 s_{\mu_i}^{-1} & \mbox{ if}\quad \mu_i>k_{i-1}.
\end{array} \right.
$$
Finally, we need to apply $\theta'$ to the above factors. Notice
that there are only two kinds of elements to consider. The first one
is $s_i$, with $i>1$, which by definition is mapped to
$\theta'(s_i)=a_{i,i-1}\, a_{i+n-1,i+n-2}$. The elements of the
second kind are those of the form $(s_{i}s_{i-1} \cdots s_2)\: s_1\:
(s_2^{-1} s_3^{-1}\cdots  s_{i}^{-1})$, for $i=1,\ldots,n-1$. One
can use the Birman-Ko-Lee presentation to show that the image under
$\theta'$ of this element is precisely $a_{i+n-1,i}$, but it is
easier to show it geometrically, since the element $(s_{i}s_{i-1}
\cdots s_2)\: s_1\: (s_2^{-1} s_3^{-1}\cdots  s_{i}^{-1})$ is
precisely the one in the right hand side of
Figure~\ref{F:conjugate_of_s_1}, in which the puncture corresponding
to the $(n-1)$-st root of unity $\xi_i$ makes a loop around the
origin. It is then easy to lift such a path via $\theta^{-1}$,
obtaining the braid $a_{i+n-1,i}$. Since $\theta^{-1}=\theta'$, one
has $\theta'\left((s_{i}s_{i-1} \cdots s_2)\: s_1\: (s_2^{-1}
s_3^{-1}\cdots  s_{i}^{-1})\right)= a_{i+n-1,i}$, as we wanted to
show.

\begin{figure}[ht]
\centerline{\includegraphics{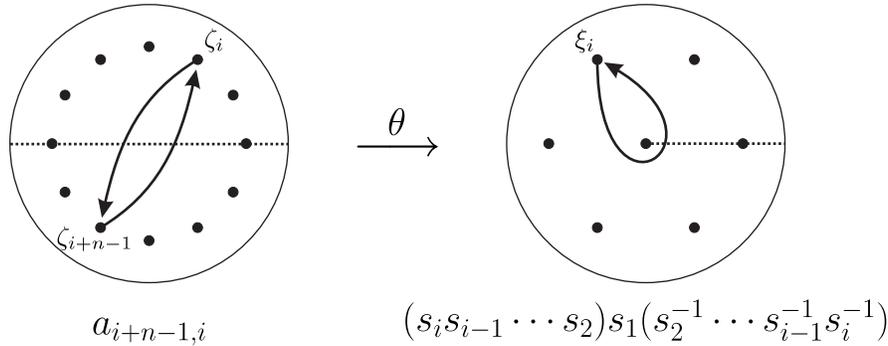}} \caption{The
image under $\theta$ of $a_{i+n-1,i}$.} \label{F:conjugate_of_s_1}
\end{figure}

One can finally transform the word
$X=\sigma_{\mu_1}^{\epsilon_1}\cdots \sigma_{\mu_l}^{\epsilon_l}$ to
a word representing $\theta'(\rho^{-1}(X))$, if one replaces each
$\sigma_{\mu_i}^{\epsilon_i}$ by $\theta'(\rho^{-1}(R_{k_{i-1}}
\sigma_{\mu_i}^{\epsilon_i} R_{k_i}^{-1}))$. By the above
discussion, the formulae in the statement hold.

It remains to notice that the numbers $\mu_i$ and $k_i$, for
$i=1,\ldots,l$ can be obtained in time $O(l)$, and that the
procedure given by the statement replaces each letter of $X$ by at
most two letters of $\theta'(\rho^{-1}(X))$. Hence the length of the
obtained word is at most $2l$, and the whole procedure has
complexity $O(l)$.
\end{proof}

Now we also need to know how to translate an element $Y\in
Sym_{2n-2}$, given as a word in the band generators of
$B_{2n-2}^B$ and their inverses, to a word representing
$\rho(\theta(Y))\in P_{n,2}\subset B_n^A$.  We first need a
preparatory result:

\begin{lemma}\label{L:from BKL to polygonal}
If $Y\in Sym_{2n-2}$ is given as a word of length $l$ in the band
generators of $B_{2n-2}^B$ and their inverses, then one can
compute in time $O(l^2n)$ a word $\delta^{t} p_1p_2\cdots p_k$
representing $Y$, such that each $p_i\in Sym_{2n-2}$ is either a
symmetric polygonal braid $\pbraid_P$, or the product of two
commuting polygonal braids $\pbraid_{P_1}\pbraid_{P_2}$ such that a
rotation of 180 degrees permutes $\pbraid_{P_1}$ and $\pbraid_{P_2}$.
Moreover, $|t|\leq l$ and $k\leq ln/2$.
\end{lemma}

\begin{proof}
The way to obtain the word $p_1\cdots p_k$ is just the computation
of the left normal form of $Y$ in $B_{2n-2}^B$. It is shown
in~\cite{Reiner} that the set of {\it symmetric} non-crossing
partitions of the $(2n-2)$-nd roots of unity (the symmetric simple
elements in $B_{2n-2}^B$) is a sublattice of the whole lattice
of non-crossing partitions. This implies that the Garside structure
of $B_{2n-2}^B$ restricts to a Garside structure on
$Sym_{2n-2}$. Therefore, since $\delta\in Sym_{2n-2}$, the greatest
common divisor of $Y$ and any power of $\delta$ is also symmetric,
and hence every factor in the left normal form of $Y$ is symmetric.

By~\cite{BKL}, the left normal form of $Y$ can be computed in time
$O(l^2n)$. Once that it is computed, each non-$\delta$ factor is the
product of mutually commuting polygonal braids, and the union of
these polygons must be symmetric. Hence, each of these polygons is
either symmetric, or it belongs of a pair of polygons which are
permuted by a rotation of 180 degrees, so the result follows.

Finally, notice that the left normal form of $Y$ has the form
$\delta^t y_1\cdots y_s$ with $|t|\leq l$ and $s\leq l$. Now every
$y_i$ contains at most one symmetric polygonal braid, namely the one
containing the origin. The remaining polygonal braids of $y_i$ come
in pairs. The symmetric polygonal braid, if it exists, involves at
least two punctures, and each pair of polygonal braids involves at
least 4 punctures. Hence $y_i$ can be decomposed into a product of
at most $1+(2n-4)/4=n/2$ factors of the form $p_j$. Since $s\leq l$,
one finally obtains $k\leq ln/2$, as we wanted to show.
\end{proof}

\begin{lemma}\label{L:from Sym to P}
Let $Y\in Sym_{2n-2}$  be given as a word of length $l$ in the band
generators and their inverses, and let $Y=\delta^t p_1\cdots p_k$ be
the decomposition given in Lemma~\ref{L:from BKL to polygonal}. Then
one obtains a word in the Artin generators and their inverses
representing $\rho(\theta(Y))$ as follows.
\begin{enumerate}

 \item Each $\delta\in B_{2n-2}^B$ should be replaced by $\rho(\theta(\delta))=\varepsilon\in B_n^A$.

\item If $p_i$ is the product of two polygonal braids
$\pbraid_{P_1}\pbraid_{P_2}$, where the vertices of the polygons are
$\{\zeta_{i_1}, \ldots, \zeta_{i_d}\}$ and $\{-\zeta_{i_1}, \ldots,
-\zeta_{i_d}\}$ respectively, let $k\in \{0,\ldots,n-2\}$ be such
that
$\{\zeta_{i_1+k},\ldots,\zeta_{i_d+k}\}=\{\zeta_{j_1},\ldots,\zeta_{j_d}\}$
with $1\leq j_1< \cdots <j_d < n$. Then $p_i$ should be replaced by
$$
\rho(\theta(\pbraid_{P_1}\pbraid_{P_2}))= \varepsilon^{k}\sigma_1
\left(\prod_{\substack{i=j_1+1 \\ (i\neq j_k\: \forall
k)}}^{j_d-1}\sigma_i^{-1} \right) (\sigma_{j_d} \sigma_{j_d-1}\cdots
\sigma_{j_1+1}) \:\sigma_1^{-1}\varepsilon^{-k}.
$$

\item If $p_i$ is a symmetric polygonal braid $\pbraid_P$, and the vertices of
the polygon $P$ are
$$
\{\zeta_{j_1}, \ldots, \zeta_{j_d},-\zeta_{j_1}, \ldots,
-\zeta_{j_d}\},
$$
with $1\leq j_1< \cdots <j_d < n$, then $p_i$ should be replaced by
$$
\rho(\theta(\pbraid_P)) = \sigma_1 \left(\prod_{\substack{i=j_1+1 \\
(i\neq j_k\: \forall k)}}^{j_d-1}\sigma_i^{-1} \right) (\sigma_{j_d}
\sigma_{j_d-1}\cdots \sigma_1) \sigma_1 (\sigma_2^{-1} \cdots
\sigma_{j_1}^{-1} )\:\sigma_1^{-1}.
$$

\end{enumerate}
\end{lemma}

\begin{proof}
Consider the element $\alpha= s_{n-1}s_{n-2}\cdots  s_1\in \mathcal
A(\mathbf B_{n-1})$. It is represented in the central picture of
Figure~\ref{F:delta_goes_to_epsilon}. On the one hand, by
Lemma~\ref{L:from Artin B to P} one has:
$$
\rho(\alpha)= \sigma_{n-1}\sigma_{n-2}\cdots \sigma_3 \, (\sigma_1
\sigma_2 \sigma_1^{-1}) \,\sigma_1^2 = \sigma_1 (\sigma_{n-1}
\sigma_{n-2} \cdots \sigma_1)= \varepsilon.
$$
On the other hand, Lemma~\ref{L:from Artin B to Sym} together with
presentation~(\ref{E:BKL presentation}) tell us that
$$
 \theta'(\alpha)=
(a_{n-1,n-2}\, a_{2n-2,2n-3})(a_{n-2,n-3}\, a_{2n-3,2n-4})\cdots
(a_{2,1}\, a_{n+1,n})\, a_{n,1}
$$
$$
 = (a_{2n-2,2n-3}\, a_{2n-3,2n-4}\cdots a_{n+1,n})(a_{n-1,n-2}\, a_{n-2,n-3}\cdots
a_{2,1})\, a_{n,1}
$$
$$
 = (a_{2n-2,2n-3}\, a_{2n-3,2n-4}\cdots a_{n+1,n})\, a_{n,n-1}\,(a_{n-1,n-2}\, a_{n-2,n-3}\cdots
a_{2,1}) = \delta.
$$
Therefore, since $\theta'=\theta^{-1}$, one has
$\rho(\theta(\delta))=\rho(\alpha)=\varepsilon$ and the first case
holds.

\begin{figure}[ht]
\centerline{\includegraphics{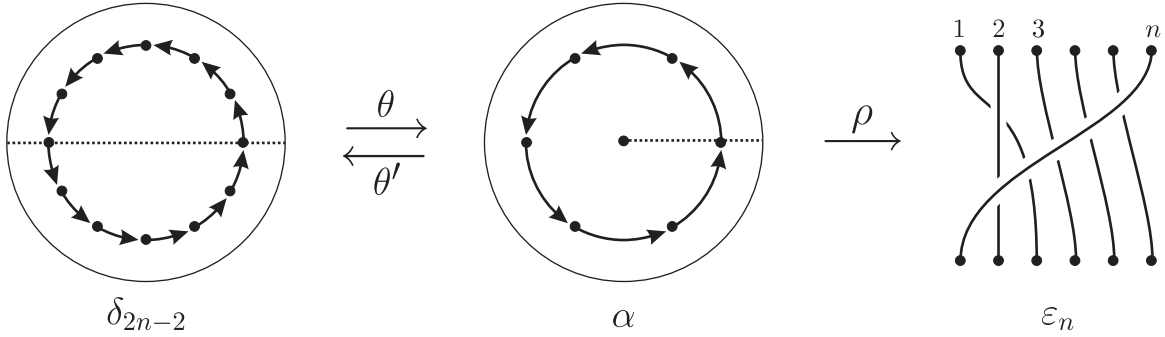}} \caption{A
geometric interpretation of $\rho(\theta(\delta))=\varepsilon$.}
\label{F:delta_goes_to_epsilon}
\end{figure}

Now suppose that $p_i$ is the product of two polygonal braids
$\pbraid_{P_1}\pbraid_{P_2}$, where the vertices of the polygons are
$\{\zeta_{i_1}, \ldots, \zeta_{i_d}\}$ and $\{-\zeta_{i_1}, \ldots,
-\zeta_{i_d}\}$. Notice that conjugation by $\delta$ in
$B_{2n-2}^B$ rotates the base points, increasing each index by
one. Therefore, since $P_1$ and $P_2$ belong to a non-crossing
partition, there exists some $k\in \{0,\ldots,n-2\}$ such that the
rotation induced by $\delta^k$ transforms $\{P_1,P_2\}$ into
$\{P_1',P_2'\}$, where the vertices of $P_1'$ belong to
$\{\zeta_1,\ldots,\zeta_{n-1}\}$.  Then
$\pbraid_{P_1}\pbraid_{P_2}=\delta^{k}
\pbraid_{P_1'}\pbraid_{P_2'}\delta^{-k}$. Since
$\rho(\theta(\delta))=\varepsilon$, in order to compute
$\rho(\theta(\pbraid_{P_1}\pbraid_{P_2}))$ it suffices to know the
value of $\rho(\theta(\pbraid_{P_1'}\pbraid_{P_2'}))$. See an example
in Figure~\ref{F:pair_of_polygons}.

\begin{figure}[ht]
\centerline{\includegraphics{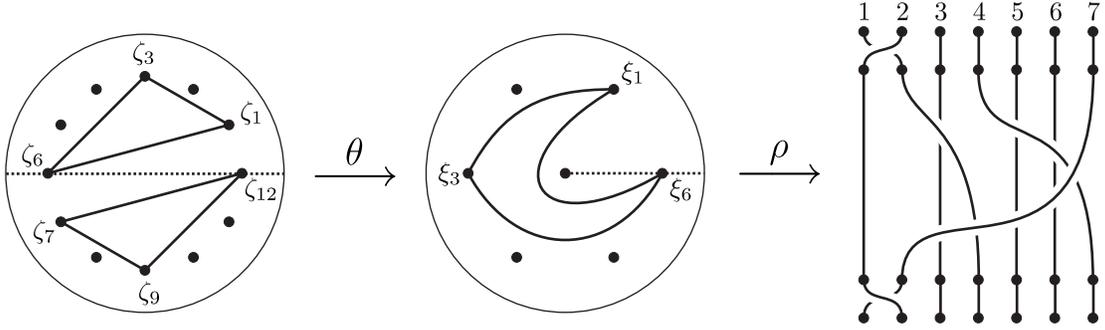}}
\caption{Translating pairs of symmetric polygonal braids in
$B_{2n-2}^B$ to Artin generators in $B_n^A$.} \label{F:pair_of_polygons}
\end{figure}

Let $\zeta_{j_1},\ldots,\zeta_{j_d}$ be the vertices of $P_1'$ in
increasing order, as in the statement. For simplicity of notation,
denote $j^*=j+n-1$ for $j=1,\ldots, n-1$. The computation goes as
follows:
$$
  \pbraid_{P_1'}\pbraid_{P_2'}= (a_{j_d,j_{d-1}}\, a_{j_{d-1},j_{d-2}}\cdots a_{j_2,j_1})
 (a_{j_d^*,j_{d-1}^*}\, a_{j_{d-1}^*,j_{d-2}^*}\cdots a_{j_2^*,j_1^*})
$$
$$
  = (a_{j_d,j_{d-1}} \, a_{j_d^*,j_{d-1}^*}) \cdots (a_{j_2,j_1} \,
a_{j_2^*,j_1^*})
$$
$$
  =\prod_{i=d}^2{(a_{j_i,j_{i-1}}\, a_{j_i^*, j_{i-1}^*})},
$$
where the index $i$ decreases from $d$ to 2.

Now one can check using Lemma~\ref{L:from Artin B to Sym} and
presentation~\ref{E:BKL presentation}, or just by drawing the
corresponding pictures, that for $1\leq u < v <n$ one has
$\theta'((s_{u+1}^{-1} s_{u+2}^{-1} \cdots s_{v-1}^{-1}) (s_{v}
s_{v-1}\cdots s_{u+1}))= a_{v,u}a_{v^*,u^*}$. Hence, since
$\theta'=\theta^{-1}$, one obtains:
$$
  \theta(\pbraid_{P_1'}\pbraid_{P_2'})=\prod_{i=d}^2{
 (s_{j_{i-1}+1}^{-1} s_{j_{i-1}+2}^{-1} \cdots s_{j_i-1}^{-1}) (s_{j_i} s_{j_i-1}\cdots
s_{j_{i-1}+1})}.
$$
Notice that $s_i$ commutes with $s_j$ if $|i-j|>1$, hence all
positive letters in the above formula can be collected to the right
(the only exception would appear if $j_{i-1}$ and $j_i$ are
consecutive for some $i$, but in that case the corresponding
negative factor is empty). It follows that:
$$
  \theta(\pbraid_{P_1'}\pbraid_{P_2'})=\left(\prod_{i=d}^2{
 (s_{j_{i-1}+1}^{-1} s_{j_{i-1}+2}^{-1}} \cdots s_{j_i-1}^{-1}) \right)
(s_{j_d} s_{j_d-1}\cdots s_{j_1+1}).
$$
Also, the $d-1$ factors made by negative letters commute with each
other, so one finally obtains:
$$
  \theta(\pbraid_{P_1'}\pbraid_{P_2'})=\left(\prod_{i=2}^d{
 (s_{j_{i-1}+1}^{-1} s_{j_{i-1}+2}^{-1}} \cdots s_{j_i-1}^{-1}) \right)
(s_{j_d} s_{j_d-1}\cdots s_{j_1+1}).
$$
$$ =
\left(\prod_{\substack{i=j_1+1 \\ (i\neq j_k\: \forall k)}}^{j_d-1}
s_i^{-1} \right) (s_{j_d} s_{j_d-1}\cdots s_{j_1+1}).
$$
Now we must apply $\rho$ to the above element. Notice that all
indices are greater than 1, so this will replace $s_2$ by $\sigma_1
\sigma_2 \sigma_1^{-1}$ and $s_i$ by $\sigma_i$ for $i>2$. This is
equivalent to replacing $s_i$ by $\sigma_1 \sigma_i \sigma_1^{-1}$ for
every $i>1$.
Hence, applying $\rho$ reduces to replacing each $s_i$ by
$\sigma_i$, and then conjugating the whole element by
$\sigma_1^{-1}$. That is,
$$
 \rho(\theta(\pbraid_{P_1'}\pbraid_{P_2'}))=
\sigma_1 \left(\prod_{\substack{i=j_1+1 \\ (i\neq j_k\: \forall
k)}}^{j_d-1}\sigma_i^{-1} \right) (\sigma_{j_d} \sigma_{j_d-1}\cdots
\sigma_{j_1+1}) \:\sigma_1^{-1},
$$
and $\rho(\theta(\pbraid_{P_1}\pbraid_{P_2}))$ is precisely as we
stated.

It remains to show the third case, in which $p_i$ is a single
symmetric polygonal braid $\pbraid_P$, where the vertices of $P$ are
$\{\zeta_{j_1},\cdots \zeta_{j_d},-\zeta_{j_1},\cdots
-\zeta_{j_d}\}= \{\zeta_{j_1},\cdots
\zeta_{j_d},\zeta_{j_1+n-1},\cdots \zeta_{j_d+n-1}\}$. An example
can be seen in Figure~\ref{F:single_polygon}.

\begin{figure}[ht]
\centerline{\includegraphics{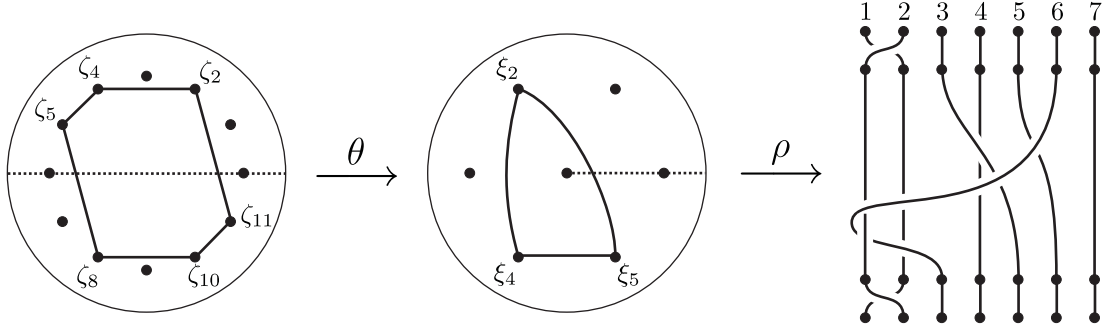}}
\caption{Translating a single symmetric polygonal braid in
$B_{2n-2}^B$ to Artin generators in $B_n^A$.} \label{F:single_polygon}
\end{figure}

Recall that $j^*=j+n-1$ for $j=1,\ldots, n-1$. In this case one has
\begin{eqnarray*}
  \pbraid_P &=& (a_{j_d^*,j_{d-1}^*}\, a_{j_{d-1}^*,j_{d-2}^*} \cdots
               a_{j_2^*,j_1^*}) \, a_{j_1^*,j_d}\, (a_{j_d,j_{d-1}}\,
               a_{j_{d-1},j_{d-2}}\cdots a_{j_2,j_1}) \\
   &=& (a_{j_d^*,j_{d-1}^*}\, a_{j_{d-1}^*,j_{d-2}^*} \cdots
       a_{j_2^*,j_1^*})  \, (a_{j_d,j_{d-1}}\, a_{j_{d-1},j_{d-2}}\cdots
       a_{j_2,j_1}) \, a_{j_1^*,j_1}.
\end{eqnarray*}

One can apply the reasoning of the previous step to the first two
factors, so it only remains to compute
$\rho(\theta(a_{j_1^*,j_1}))$. This is done by noticing that
\begin{eqnarray*}
  a_{j_1^*,j_1} &=&
        (a_{j_1,j_1-1}\,a_{j_1^*,j_1^*-1})
        (a_{j_1-1,j_1-2}\,a_{j_1^*-1,j_1^*-2})\cdots
        (a_{2,1}\,a_{n+1,n}) \cdot
        a_{n,1} \cdot \\
     & & \cdot
        (a_{2,1}^{-1}\,a_{n+1,n}^{-1})\cdots
        (a_{j_1-1,j_1-2}^{-1}\,a_{j_1^*-1,j_1^*-2}^{-1})
        (a_{j_1,j_1-1}^{-1}\,a_{j_1^*,j_1^*-1}^{-1})
                  ,
\end{eqnarray*}
which yields
$$
  \theta(a_{j_1,j_1^*})= (\theta')^{-1}(a_{j_1,j_1^*})=
  (s_{j_1}\cdots s_2) \, s_1 \, (s_2^{-1}\cdots s_{j_1}^{-1}).
$$
Since applying $\rho$ reduces to replacing $s_1$ by $\sigma_1^2$,
then $s_i$ by $\sigma_i$ for $i>1$, and then conjugating everything
by $\sigma_1^{-1}$, one obtains:
$$
  \rho(\theta(a_{j_1,j_1^*}))=
  \sigma_1 (\sigma_{j_1} \cdots \sigma_2) \, \sigma_1^2 \,
  (\sigma_2^{-1}\cdots \sigma_{j_1}^{-1}) \sigma_1^{-1}.
$$

Therefore
$$
  \rho(\theta(\pbraid_P))= \sigma_1 \left(\prod_{\substack{i=j_1+1 \\ (i\neq j_k\: \forall
k)}}^{j_d-1}\sigma_i^{-1} \right) (\sigma_{j_d} \sigma_{j_d-1}\cdots
\sigma_{j_1+1}) (\sigma_{j_1} \cdots \sigma_2) \, \sigma_1^2 \,
(\sigma_2^{-1}\cdots \sigma_{j_1}^{-1}) \sigma_1^{-1},
$$
which is precisely the formula in the statement, so the proof is
finished.
\end{proof}

\subsubsection{Using symmetric braids to solve the conjugacy search problem}
\label{SS-3.1.2}

Recall that we are given $X\in B_n^A$ as a word in the Artin
generators $\sigma_1,\ldots, \sigma_{n-1}$ and their inverses, and
we know that $X$ is conjugate to $\varepsilon^k$ for some $k\neq 0$.
This means that the permutation $\pi_X$ consists of the $k$-th power
of a cycle of length $n-1$, that is $\pi_X=(a)(b_1\ \cdots\
b_{n-1})^k$, where $a\neq b_i$ for every $i$.

The easy case happens when $k$ is a multiple of $n-1$, say
$k=(n-1)t$. Then $\varepsilon^k=\Delta^{2t}$, so $X$ is conjugate to
a power of $\Delta^2$. But since $\Delta^2$ is a central element,
this implies that $X=\Delta^{2t}$. Hence $X=\varepsilon^k$ and we
are done.

We can then assume that $k$ is not a multiple of $n-1$. This means
that the only puncture which is fixed by $X$ is the $a$-th one. If
we denote $C_1=\sigma_{[a\rightarrow 2]}$, it clearly follows that
$Y=C_1^{-1}XC_1$ fixes the second strand, that is, $Y\in P_{n,2}$.
Notice also that $\varepsilon\in P_{n,2}$, so $\varepsilon^k\in
P_{n,2}$. This means that $Y$ and $\varepsilon^k$ are two elements
in $P_{n,2}$ which are conjugate in $B_n^A$. Fortunately, they are
also conjugate in $P_{n,2}$, as it is shown in the following result.

\begin{lemma}
If $Y,Z\in P_{n,2}$ are conjugate braids whose permutations have a
single fixed point (namely 2), then for every conjugating element
$C\in B_n^A$ such that $C^{-1}YC=Z$, one has $C\in P_{n,2}$.
\end{lemma}

\begin{proof}
Let $j=\pi_C(2)$. If $j\neq 2$, then $\pi_{YC}(2) =
\pi_{C}(\pi_{Y}(2)) = \pi_C(2)=j$, while $\pi_{CZ}(2) =
\pi_Z(\pi_C(2)) = \pi_Z(j)\neq j$ (since the only fixed point of
$\pi_Z$ is 2, and $j\neq 2$). This contradicts the
assumption $YC=CZ$, so we must have $\pi_C(2)=2$, that is $C\in
P_{n,2}$.
\end{proof}

As a consequence, every conjugating element from $Y$ to
$\varepsilon^k$, when $k$ is not a multiple of $n-1$, must belong to
$P_{n,2}$. Therefore, finding a conjugating element from $Y$ to
$\varepsilon^k$ in $B_n^A$ reduces to solving the conjugacy search
problem in $P_{n,2}$ for conjugates of $\varepsilon^k$.

Our strategy consists of applying $\theta' \circ \rho^{-1}$, solving
the resulting problem in $Sym_{2n-2}$, and then mapping the solution
back to $P_{n,2}$ using $\rho\circ \theta$. Recall from
Lemma~\ref{L:from Sym to P} that $\rho(\theta(\delta))=\varepsilon$,
hence $\theta' (\rho^{-1} (\varepsilon))= \delta \in Sym_{2n-2}$.
Therefore we must solve the conjugacy search problem in $Sym_{2n-2}$
for $\theta'( \rho^{-1}(Y))$ and $\delta^k$.

Recall that, as a consequence of~\cite{Reiner}, the group
$Sym_{2n-2}$ has a Garside structure which is the restriction of the
Birman-Ko-Lee structure of $B_{2n-2}^B$. The Garside element
of this structure is hence $\delta$, so the conjugacy search problem
for powers of $\delta\in Sym_{2n-2}$ can be solved very fast, by
applying iterated cyclings and decyclings.  But one does not need to
care about the Garside structure of $Sym_{2n-2}$, since one can
directly work with the Garside structure of $B_{2n-2}^B$, as
it is shown in the following result.

\begin{lemma}\label{L:cycling and decycling symmetric braids}
Let $Z\in Sym_{2n-2} \subset B_{2n-2}^B$ be given as a word of
length $l$ in the band generators and their inverses. Suppose that
$Z$ is conjugate to $\delta^k$ for some $k\neq 0$. Then by applying
at most $(2n-3)l$ cyclings and decyclings to $Z$, using the Garside
structure of $B_{2n-2}^B$, one conjugates $Z$ to $\delta^k$
and the conjugating element that is obtained belongs to $Sym_{2n-2}$.
\end{lemma}

\begin{proof}
By~\cite{BKL-2001}, by applying at most $(2n-3)l$ cyclings and
decyclings to $Z$ one obtains an element which has minimal canonical
length. Since $Z$ is conjugate to $\delta^k$, and $\delta$ is the
Garside element of $B_{2n-2}^B$, it follows that the resulting
element is precisely $\delta^k$. Hence one obtains $C\in
B_{2n-2}^B$ such that $C^{-1}ZC=\delta^k$.

Now recall that if a braid in $B_{2n-2}^B$ is symmetric, then
every factor in its left normal form is also symmetric. Hence the
conjugating elements in all cyclings and decyclings applied above
are symmetric braids, so $C\in Sym_{2n-2}$, as we wanted to show.
\end{proof}

This finally gives us the algorithm to solve the conjugacy search
problem for conjugates of $\varepsilon^k$.

\medskip
\noindent {\bf Algorithm C:}

Input: A word $w$ in Artin generators and their inverses
representing $X\in B_n^A$ conjugate to $\varepsilon^k$.

Output: $C\in B_n^A$ such that $C^{-1}XC=\varepsilon^k$.

\begin{enumerate}

 \item If $k$ is a multiple of $n-1$, return $C=1$.

 \item Compute $a$, the only puncture fixed by $\pi_X$. Let
$Y=\sigma_{[a\rightarrow 2]}^{-1} X \sigma_{[a\rightarrow 2]}\in
P_{n,2}$.

 \item Using Lemma~\ref{L:from P to Sym}, compute
$Z=\theta'(\rho^{-1}(Y))$.

 \item Apply iterated cycling and decycling to $Z\in B_{2n-2}^B$ until $\delta^k$
is obtained. Let $C_0\in Sym_{2n-2}$ be the conjugating element.

 \item Using Lemma~\ref{L:from Sym to P}, compute
$C_1=\rho(\theta(C_0))$.

 \item Return $C=\sigma_{[a\rightarrow 2]}C_1$.

\end{enumerate}

\medskip

\begin{proposition}\label{P:Algorithm C}
 Algorithm C has complexity $O(l^3n^2)$.
\end{proposition}

\begin{proof}
 The number $a$ in step 2 can be computed in time $O(ln)$, and the
letter length of the word $Y$ is at most $2(n-2)+l$, hence, by
Lemma~\ref{L:from P to Sym}, the word $Z$ is obtained in time
$O(n+l)$, and its letter length is at most $4(n-2)+2l$, that is,
$O(n+l)$ as well. By Lemma~\ref{L:n<l}, $O(n+l)=O(l)$, so the letter
length of $Z$ in band generators is $O(l)$.

By Lemma~\ref{L:cycling and decycling symmetric braids}, one just
needs to apply $O(nl)$ cyclings and decyclings to $Z$ in step 4,
each computation taking time $O(l^2n)$ since it is equivalent to
computing a left normal form of a word of length $O(l)$. Hence, step
4 takes time $O(l^3n^2)$, and it is the most time-consuming step of
the algorithm. The conjugating element $C_0\in Sym_{2n-2}$ consists
of at most $O(ln)$ simple factors.

Notice that $C_0$ is already given as a product of symmetric simple
elements. Hence one can directly apply the formulae in
Lemma~\ref{L:from Sym to P}, to compute $C_1=\rho(\theta(C_0))$.
Since there are $O(ln)$ factors, and each one is replaced by at most
$n/2$ words of letter length  $O(3n+2n^2)=O(n^2)$, it follows that
step 5 takes time $O(ln^4)=O(l^3n^2)$, hence the whole algorithm has
complexity $O(l^3n^2)$ as we wanted to show.
\end{proof}

\subsection{The complete algorithm}

We are finally ready to prove Theorem~\ref{T:main theorem} by giving an
algorithm which solves step (3) in the statement of
Theorem~\ref{T:main theorem} in time $O(l^3n^2\log n)$.


\medskip
\noindent {\bf Algorithm D}

Input: Two words $w_X$, $w_Y$ in Artin generators and their inverses
representing two braids $X,Y \in B_n^A$.

Output: `Fail' if either $X$ or $Y$ is not periodic, or if they are
not conjugate. Otherwise, an element $C\in B_n^A$ such that $C^{-1} X C = Y$.

\begin{enumerate}

 \item Apply Algorithm A to $w_X$ and $w_Y$.

 \item If either $X$ or $Y$ is not periodic return `Fail'. If $X$
and $Y$ are not conjugate to the same power of $\delta$ or
$\varepsilon$, return `Fail'.

 \item If $X$ and $Y$ are conjugate to $\delta^k$ for some $k$,
apply Algorithm $B$ to $X$ and $Y$ to find $C_1,C_2\in B_n^A$ such
that $C_1^{-1}XC_1=\delta^k=C_2^{-1}YC_2$. Return $C=C_1C_2^{-1}$.

 \item If $X$ and $Y$ are conjugate to $\varepsilon^k$ for some $k$,
apply Algorithm $C$ to $X$ and $Y$ to find $C_1,C_2\in B_n^A$ such
that $C_1^{-1}XC_1=\varepsilon^k=C_2^{-1}YC_2$. Return
$C=C_1C_2^{-1}$.

\end{enumerate}

\medskip

\begin{proposition}
Algorithm D has complexity $O(l^3n^2\log n)$, where
$l=\max\{|w_X|,|w_Y|\}$.
\end{proposition}

\begin{proof}
By Proposition~\ref{P:Algorithm A}, Algorithm A has complexity
$O(l^2n^3\log n)$. By Proposition~\ref{P:Algorithm B}, the
complexity of Algorithm B is $O(l^3n^2)$, which is the same
complexity as that of Algorithm C, by Proposition~\ref{P:Algorithm
C}. Therefore, Algorithm D has complexity $O(l^2n^3\log n +
l^3n^2)$. By Lemma~\ref{L:n<l}, this complexity is equivalent to
$O(l^3n^2\log n)$, as we wanted to show.
\end{proof}

\section{Timing results}
\label{S:Timing results}
In this section we present and analyze running times for the conjugacy
search for periodic elements in Artin braid groups;  we compare the
established algorithm based on computing ultra summit
sets~\cite{Gebhardt} to the algorithms developed in this paper.

For several values of the parameters $n$, $k$ and $c$, tests in $B_n$ were
conducted as follows.
\begin{enumerate}
\item For $i=1,\dots,100$, we construct a pseudo-random element $z_i\in B_n^A$
  as the product of $c$ randomly chosen simple elements.
\item We compute the samples $\{(\delta^k)^{z_i} : i=1,\dots,100\}$ and
  $\{(\varepsilon^k)^{z_i} : i=1,\dots,100\}$; each element is stored in left
  normal form.
\item \label{step:computation} For each element $x$ in a sample we compute an
  element conjugating $x$ to $\delta^k$ or $\varepsilon^k$, respectively.
\end{enumerate}

Step \ref{step:computation} was performed separately for each sample,
first using the algorithm from \cite{Gebhardt}, in the sequel referred
to as Algorithm U, and then again using Algorithm B or Algorithm C.
Only the total time for this step was measured for each case.  A
memory limit of 512\,MB and a time limit of 250\,s were applied for
each test.

All computations were performed on a Linux PC with a 2.4\,GHz Pentium
4 CPU, 533\,MHz system bus and 1.5\,GB of RAM using the computational
algebra system \Magma{}~\cite{magma}.   An implementation of Algorithm
U written in C is part of the \Magma{} kernel; Algorithms B and C were
implemented in the \Magma{} language.

\noindent {\bf Remark: } One technical aspect of the
implementation of Algorithms B and C needs to be mentioned briefly to
explain the observed behavior.

As Algorithms B and C involve mapping a given word, generator by
generator, to another Garside group, a naive implementation of these
algorithms will react very sensitively to the word length of the given
element $x$.

Note, however, that a conjugate $y$ of $x$ having minimal canonical
length with respect to the usual Garside structure, together with
a conjugating element, can be computed by iterated application of
cycling and decycling in time $O(\ell^3n^3\log n)$, where $\ell$ is
the number of simple factors of $x$.\footnote{Note that $\ell$, unlike
  the letter length $l$, is not bounded below by $n$ for periodic
  braids.}
Note further that if $x$ is periodic, the canonical length of $y$ as
above is at most 1.
Moreover, powers
of $\Delta^2$ can be discarded for the purpose of computing
conjugating elements, as $\Delta^2$ is central in $B_n$.  The
techniques from Algorithms B and C are then applied to the resulting
element whose length in terms of Artin generators is bounded by $n^2$.

While this does not improve the complexity bounds, it significantly
reduces computation times, especially for large values of the
parameter $c$ above, and is critical for the cross-over points between
Algorithm U on the one hand and Algorithms B and C on the other hand.

We finally remark that in the special case that the minimal canonical
length of conjugates of $x$ is 0, that is, in the case that $x$ is
conjugate to a power of $\Delta$, its ultra summit set has cardinality
1 and we do not have to use Algorithms B and C, as a conjugating
element can be obtained directly, just by iterated application of
cycling and decycling.
\bigskip

\begin{table}
 \caption{Total execution times of Algorithms U, B and C for all 100
   elements of a sample for $c=10$ and various values of $n$ and $k$.
   Where no value is given, either the memory limit of 512\,MB or the
   time limit of 250\,s was exceeded.}
 \label{table_c10}
 \medskip
 \begin{tabular}{c||
       c|c|c|c|c|c||
       c|c|c|c|c|c}
  $k$
       & \multicolumn{6}{c||}{1}
       & \multicolumn{6}{c}{2} \\
  \hline
  $n$
      & 5 & 7 & 10 & 15 & 20 & 50
      & 5 & 7 & 10 & 15 & 20 & 50 \\
  \hline
  U[$\delta$]
      & 0.03 & 0.12 & 1.56 & 88.14 & --- & ---
      & 0.02 & 0.38 & 22.15 & --- & --- & ---
  \\
  B
      & 0.02 & 0.04 & 0.07 & 0.16 & 0.34 & 3.56
      & 0.02 & 0.03 & 0.06 & 0.16 & 0.29 & 2.75
  \\
  U[$\varepsilon$]
      & 0.03 & 0.19 & 4.05 & --- & --- & ---
      & 0.02 & 0.16 & 64.22 & --- & --- & ---
  \\
  C
      & 0.05 & 0.12 & 0.23 & 0.53 & 0.97 & 6.92
      & 0.01 & 0.10 & 0.25 & 0.52 & 1.01 & 6.95
  \\
 \end{tabular}
 \smallskip
 \begin{tabular}{c||
       c|c|c|c|c||
       c|c|c|c||
       c|c|c}
  $k$
       & \multicolumn{5}{c||}{3}
       & \multicolumn{4}{c||}{4}
       & \multicolumn{3}{c}{6} \\
  \hline
  $n$
      & 7 & 10 & 15 & 20 & 50
          & 10 & 15 & 20 & 50
               & 15 & 20 & 50 \\
  \hline
  U[$\delta$]
      & 0.05 & 58.81 & --- & --- & ---
             & 3.86 & --- & --- & ---
                    & --- & --- & ---
  \\
  B
      & 0.04 & 0.08 & 0.12 & 0.34 & 2.79
             & 0.06 & 0.16 & 0.23 & 2.37
                    & 0.10 & 0.29 & 2.39
  \\
  U[$\varepsilon$]
      & 0.02 & 9.59 & --- & --- & ---
             & 0.45 & --- & --- & ---
                    & --- & --- & ---
  \\
  C
      & 0.01 & 0.22 & 0.60 & 1.03 & 7.02
             & 0.22 & 0.57 & 1.07 & 7.02
                    & 0.53 & 1.09 & 7.22
  \\
 \end{tabular}
 \smallskip
 \begin{tabular}{c||
       c|c|c||
       c|c||
       c|c||
       c|c||
       c||
       c}
  $k$
       & \multicolumn{3}{c||}{7}
       & \multicolumn{2}{c||}{8}
       & \multicolumn{2}{c||}{9}
       & \multicolumn{2}{c||}{10}
       & \multicolumn{1}{c||}{11}
       & \multicolumn{1}{c}{12} \\
  \hline
  $n$
      & 15 & 20 & 50
           & 20 & 50
           & 20 & 50
           & 20 & 50
                & 50
                & 50 \\
  \hline
  U[$\delta$]
      & 6.17 & ---  & ---
             & ---  & ---
             & ---  & ---
             & 0.16 & ---
                    & ---
                    & ---
  \\
  B
      & 0.12 & 0.33 & 3.04
             & 0.18 & 2.60
             & 0.23 & 3.03
             & 0.03 & 1.71
                    & 3.18
                    & 2.83
  \\
  U[$\varepsilon$]
      & 0.09 & --- & ---
             & --- & ---
             & 130.34 & ---
             & 67.69  & ---
                    & ---
                    & ---
  \\
  C
      & 0.02 & 1.06 & 7.86
             & 1.02 & 7.68
             & 0.95 & 7.84
             & 0.73 & 7.96
                    & 8.23
                    & 8.26
 \end{tabular}
\end{table}

\begin{table}
 \caption{Total execution times of Algorithms U, B and C for all 100
   elements of a sample for $c=250$ and various values of $n$ and $k$.
   Where no value is given, either the memory limit of 512\,MB or the
   time limit of 250\,s was exceeded.}
 \label{table_c250}
 \medskip
 \begin{tabular}{c||
       c|c|c|c|c|c||
       c|c|c|c|c|c}
  $k$
       & \multicolumn{6}{c||}{1}
       & \multicolumn{6}{c}{2} \\
  \hline
  $n$
      & 5 & 7 & 10 & 15 & 20 & 50
      & 5 & 7 & 10 & 15 & 20 & 50 \\
  \hline
  U[$\delta$]
      & 0.16 & 0.40 & 2.05 & 85.20 & --- & ---
      & 0.15 & 0.65 & 20.42 & --- & --- & ---
  \\
  B
      & 0.16 & 0.32 & 0.67 & 1.21 & 1.83 & 8.24
      & 0.16 & 0.33 & 0.66 & 1.22 & 1.76 & 6.79
  \\
  U[$\varepsilon$]
      & 0.16 & 0.49 & 4.32 & --- & --- & ---
      & 0.14 & 0.42 & 59.76 & --- & --- & ---
  \\
  C
      & 0.19 & 0.40 & 0.83 & 1.51 & 2.37 & 10.75
      & 0.14 & 0.38 & 0.86 & 1.57 & 2.42 & 10.69
  \\
 \end{tabular}
 \smallskip
 \begin{tabular}{c||
       c|c|c|c|c||
       c|c|c|c||
       c|c|c}
  $k$
       & \multicolumn{5}{c||}{3}
       & \multicolumn{4}{c||}{4}
       & \multicolumn{3}{c}{6} \\
  \hline
  $n$
      & 7 & 10 & 15 & 20 & 50
          & 10 & 15 & 20 & 50
               & 15 & 20 & 50 \\
  \hline
  U[$\delta$]
      & 0.33 & 56.14 & --- & --- & ---
             & 4.36 & --- & --- & ---
                    & --- & --- & ---
  \\
  B
      & 0.31 & 0.69 & 1.14 & 1.81 & 6.86
             & 0.65 & 1.22 & 1.66 & 6.26
                    & 1.11 & 1.76 & 6.36
  \\
  U[$\varepsilon$]
      & 0.31 & 9.64 & --- & --- & ---
             & 0.99 & --- & --- & ---
                    & --- & --- & ---
  \\
  C
      & 0.29 & 0.83 & 1.59 & 2.47 & 11.06
             & 0.85 & 1.60 & 2.55 & 10.85
                    & 1.57 & 2.52 & 11.19
  \\
 \end{tabular}
 \smallskip
 \begin{tabular}{c||
       c|c|c||
       c|c||
       c|c||
       c|c||
       c||
       c}
  $k$
       & \multicolumn{3}{c||}{7}
       & \multicolumn{2}{c||}{8}
       & \multicolumn{2}{c||}{9}
       & \multicolumn{2}{c||}{10}
       & \multicolumn{1}{c||}{11}
       & \multicolumn{1}{c}{12} \\
  \hline
  $n$
      & 15 & 20 & 50
           & 20 & 50
           & 20 & 50
           & 20 & 50
                & 50
                & 50 \\
  \hline
  U[$\delta$]
      & 7.72 & --- & ---
             & --- & ---
             & --- & ---
             & 1.44 & ---
                    & ---
                    & ---
  \\
  B
      & 1.15 & 1.89 & 6.79
             & 1.62 & 6.37
             & 1.70 & 6.80
             & 1.41 & 5.23
                    & 7.00
                    & 6.53
  \\
  U[$\varepsilon$]
      & 1.04 & --- & ---
             & --- & ---
             & 162.83 & ---
             & 90.88 & ---
                    & ---
                    & ---
  \\
  C
      & 0.99 & 2.50 & 11.57
             & 2.49 & 11.47
             & 2.43 & 11.55
             & 2.21 & 11.70
                    & 12.51
                    & 11.93
 \end{tabular}
\end{table}

The main results can be summarized as follows; see
Tables~\ref{table_c10} and \ref{table_c250}.
\begin{enumerate}
\item Time (and memory) requirements of Algorithm U increase rapidly
  with increasing value of $n$.  With the exception of elements which
  are conjugate to a power of $\Delta$, conjugacy search using
  Algorithm U fails for $n\gtrsim 15$.

  In the light of the exponential growth of $USS(\delta)$ and
  $USS(\varepsilon)$ established in Corollaries~\ref{C:size of
  USS(delta)} and \ref{C:size of USS(epsilon)}, this had to be
  expected.
\item In contrast to this, the computation times for Algorithms B and
  C grow much more slowly with increasing value of $n$.  The data is
  consistent with a polynomial growth; a regression analysis for fixed
  values of the parameters $k$ and $c$ suggests that average times are
  proportional to $n^{e_n}$, where the value $e_n\approx 2.2$ is suggested by a regression analysis.\footnote{Note that
  for fixed values of $k$ and $c$ the word length $l$ is not fixed but
  grows at least linearly in $n$; cf.~Lemma~\ref{L:n<l}.  Hence this
  value of $e_n$ does not contradict the complexity bounds from
  Propositions~\ref{P:Algorithm B} and \ref{P:Algorithm C}.}

  In particular, solving the conjugacy search problem for periodic
  elements using Algorithm D is is not significantly harder than other
  operations in with braids, that is, it is feasible whenever the parameter
  values permit any computations at all.
\item The computation times of Algorithm U depend in a very sensitive
  way on the value of $k$, whereas the running times of Algorithms B
  and C, with the exception of elements which are conjugate to a power
  of $\Delta$ and are treated differently, show relatively little
  dependency on $k$.
\item Average running times for all algorithms appear to be sub-linear
  in $c$ for fixed values of the parameters $n$ and $k$.

  For Algorithm U, the effect of $c$ becomes negligible for $n\gtrsim
  10$.  This is no surprise as the value of $c$ only affects the
  initial computation of a conjugate with minimal canonical length;
  the time used in this step of the computation is only relevant if
  the ultra summit set is small.
\item Using the implementations as explained above, the cross-over
  point between Algorithm U and Algorithm B was $n\approx 5$, whereas
  the cross-over point between Algorithm U and Algorithm C was
  $n\approx 7$; the latter corresponds to the cross-over point
  between Algorithm U and Algorithm D for the implementations used in
  our tests.

  We remark that the fact that Algorithms B and C were implemented in
  the \Magma{} language (which is partly an interpreter language)
  incurs some overhead compared to the C implementation of
  Algorithm U.  This overhead is probably not significant for
  Algorithm B, as its implementation is quite simple.\footnote{Uses 20
  lines of \Magma{} code. As \Magma{} provides a kernel function
  computing ultra summit sets with respect to the Birman-Ko-Lee
  presentation, no low level operations had to be written in the
  \Magma{} language.}
  However, for Algorithm C the overhead can be expected to be
  significant, as its implementation is rather
  complex.\footnote{Uses 200 lines of \Magma{} code. Many low level
  operations had to be written in the \Magma{} language.}   This
  difference can be assumed to be the main cause for the different
  cross-over points, whence a cross-over point of $n\approx 5$ for
  comparable implementations of Algorithms U and D seems likely.
\end{enumerate}

\begin{remark}
{\rm  After this paper was accepted for publication, and as we were preparing this final copy for the publisher, we learned that E-K Lee and S.J. Lee  had posted on the arXiv their own solution to the same problem \cite{L-L}.  They reference our work and suggest some small improvements in it.  }
\end{remark}

\footnotesize
\hspace*{-8mm}
\begin{tabular}{lll}
 {\bf Joan S. Birman}  &  {\bf Volker Gebhardt}  &  {\bf Juan Gonz\'alez-Meneses} \\
 Department of Mathematics, & School of Computing and Mathematics, & Departamento de \'Algebra, \\
 Barnard College and Columbia University, & University of Western Sydney, & Universidad de Sevilla,\\
 2990 Broadway, & Locked Bag 1797, &  Apdo. 1160, \\
 New York, New York 10027, USA. & Penrith South DC NSW 1797, Australia, & 41080 Sevilla, Spain.\\
 {\tt jb@math.columbia.edu} & {\tt v.gebhardt@uws.edu.au} & {\tt meneses@us.es}
\end{tabular}

\end{document}